\documentclass[11pt]{article}
\usepackage{multirow}
\usepackage{amsmath}
\usepackage{amsopn}
\usepackage{amsthm}
\usepackage{amsfonts}
\usepackage[T1]{fontenc}
\usepackage[latin1]{inputenc}
\usepackage{graphicx,color}
\usepackage{epsfig}
\DeclareGraphicsExtensions{.eps,.ps}
\usepackage{amssymb}
\usepackage[frenchb,english]{babel}
\usepackage{url}
\usepackage[mathscr]{eucal}
\usepackage{amsbsy}
\usepackage{fancyhdr}
\usepackage{hyperref}
\usepackage{titlesec}

\numberwithin{equation}{section}

\newcounter{ex}
\newcommand\exemple[1]{\refstepcounter{ex}\par\rm{\bf Exemple \arabic{ex}.}\quad\label{#1}}

\newcounter{eq}

\newcounter{rq}
\newcommand\remarque{\refstepcounter{rq}\par\rm{\bf Remarque \arabic{rq}.}\quad}

\theoremstyle{plain}
\newtheorem{theoreme}{Th\'{e}or\`{e}me}[section]

\newtheorem{definition}[theoreme]{D\'{e}finition}

\def\one{\hbox{1\hskip -3pt I}}
\def\rd{\Bbb{R}^d}
\def\sd{S^{d-1}}

\def\sas{{$\mathcal{S}t\alpha\mathcal{S}$}}
\def\cok{\hbox{I\hskip -2pt K}}
\def\reip{{\bf P}}

\def\reixs{{\bf x}}
\def\reiys{{\bf y}}

\def\reia{{\bf a}}
\def\reib{{\bf b}}

\def\reie{{\bf E}}

\def\cvl1{\stackrel{\mathcal{L}_1}{\longrightarrow}}

\def\ints{\int_{S^{d-1}}}
\def\egenloi{\stackrel{\mathcal{L}}{=}}

\title{Illustration par quelques exemples des lois strictement stables dans un cône convexe}
\author{%
\renewcommand{\thefootnote}{\alph{footnote}}
Shuyan \textsc{Liu}\,\footnotemark[1]{}
}
\date{}
\begin{document}
\maketitle
\renewcommand{\thefootnote}{\alph{footnote}\,}
\footnotetext[1]{SAMM (EA 4543),
Université Paris 1 Panthéon-Sorbonne, Centre PMF,
75013 Paris, France.}

\renewcommand{\thefootnote}{\arabic{footnote}}

\selectlanguage{frenchb}
\begin{abstract}
La stabilité d'une variable aléatoire peut être généralisée dans un cône convexe. Les résultats concernant la représentation de LePage et le domaine d'attraction sont analogues mais différents de ceux qui sont obtenus dans l'espaces de Banach. Quelques exemples de lois strictement max-stables et $\alpha$-stables sont présentés dans cet article afin de montrer le lien entre le processus ponctuel poissonnien et la loi stable dans un cône convexe.
\end{abstract}


\smallskip\noindent \textit{Mots clés :} $\alpha$-stable, max-stable, cône convexe, processus ponctuel poissonnien, série de LePage, simulation.



\section{Introduction}\label{introduction}

Un vecteur aléatoire $X$ à valeurs dans un espace de Banach a la loi {\em strictement $\alpha$-stable} (\sas) si pour tous $a,b > 0$
\begin{equation}\label{stabintro}
a^{1/\alpha}X_1+b^{1/\alpha}X_2\egenloi (a+b)^{1/\alpha}X, 
\end{equation}
où $X_1, X_2$ sont les copies indépendantes de $X$, et $\egenloi$ représente l'égalité en loi. La stabilité (\ref{stabintro}) peut être définie dans tous les espaces où l'addition des éléments et la multiplication par les scalaires positifs sont bien définies, i.e. dans un cône convexe. Dans ce cas les résultats principaux (la représentation de LePage, le domaine d'attraction, etc.) sont analogues mais différents de ceux qui sont obtenus dans l'espace de Banach. L'intervalle des valeurs possibles de l'indice caractéristique $\alpha$ dépend des propriétés algébriques et topologiques de base de l'espace où l'élément aléatoire \sas \, est défini. Par exemple il est connu que $\alpha$ est dans l'intervalle $(0, 2]$ pour les lois $\alpha$-stables dans l'espace de Banach. Si on remplace l'addition dans (\ref{stabintro}) par les opérations de maximum et minimum, les vecteurs aléatoires max-stables et min-stables apparaissent. Les lois max-stables existent pour tout $\alpha >0$, tandis que les lois min-stables peuvent avoir toutes les valeurs négatives pour l'indice $\alpha$. Davydov et ses co-auteurs ont caractérisé dans \cite{Davydov08} les lois des éléments aléatoires \sas s dans un cône convexe et présenté les valeurs possibles pour l'indice $\alpha$ selon les propriétés du semigroupe et la distributivité de l'opération de multiplication correspondante. Nous présentons dans cet article quelques exemples des lois \sas s et max-stables dans $\Bbb{R}^2$ avec la simulation numérique de la densité, afin d'illustrer les rôles de l'indice $\alpha$ et de la mesure spectrale dans la forme de la densité de loi \sas \, dans différents cônes convexes. 

Rappelons qu'un vecteur aléatoire $\alpha$-stable symétrique ou avec $\alpha<1$ à valeurs dans $\rd$ admet la représentation de LePage suivante \cite{Lepage81}
\begin{equation}\label{lepageintro}
\sum_{k=1}^\infty c\Gamma_k^{-1/\alpha}\epsilon_k,
\end{equation}
où $c>0$, $\Gamma_k$ représente la somme de $k$ variables i.i.d. de loi exponentielle d'espérance $1$ et $\epsilon_1,\epsilon_2,\ldots$ sont les vecteurs aléatoires unitaires i.i.d. et indépendants de $\{\Gamma_k\}_{k=1}^\infty$. Notons $\Pi_\alpha$ le processus ponctuel dont le support est $\{c\Gamma_k^{-1/\alpha}\epsilon_k,k\geq 1\}$, i.e.
\begin{equation}\label{ppintro}
\Pi_{\alpha}= \sum_{k=1}^{\infty}\delta_{c\Gamma_k^{-1/\alpha}\epsilon_k}
\end{equation}
où $\delta_x$ désigne la mesure de Dirac au point $x$. Le support du processus $\Pi_\alpha$ est un ensemble fermé aléatoire union-stable \cite{Molchanov93} \cite{Davydov08}, i.e. pour tous $a,b>0$
\begin{equation}\label{stabpi}
a^{1/\alpha}\mbox{supp}(\Pi_{\alpha}')\cup b^{1/\alpha}\mbox{supp}(\Pi_{\alpha}'')\egenloi (a+b)^{1/\alpha}\mbox{supp}(\Pi_{\alpha}).
\end{equation}
où $\Pi_\alpha'$ et $\Pi_\alpha''$ sont les copies indépendantes de $\Pi_\alpha$. Cela signifie que l'union-stabilité du support de processus $\Pi_\alpha$ est étroitement liée à la stabilité classique (\ref{stabintro}). Cette relation est l'idée clé dans \cite{Davydov08} pour obtenir les résultats concernant la représentation de LePage et le domaine d'attraction de loi \sas \, dans un cône convexe. Les exemples présentés dans cet article montrent que les processus ponctuels sont utiles non seulement pour obtenir les résultats généraux sur les lois \sas s dans un cône convexe, mais aussi pour déduire la fonction de répartition et la régularité de la queue d'une loi max-stable avec la mesure spectrale donnée.

Un autre but de cet article est de montrer que les processus ponctuels poissonniens peuvent être utilisés pour simuler les vecteurs aléatoires suivant les lois \sas s et max-stables. Ces deux familles de lois apparaissent souvent dans les applications du domaine des événements rares qui demandent le traitement de valeurs extrêmes. Il existe des méthodes connues pour engendrer les données de loi stable \cite{Chambers76} \cite{Weron96} \cite{Modarres94} et de loi max-stable \cite{Nadarajah99} \cite{Stephenson03} \cite{Shi92}. En utilisant la représentation de LePage généralisée dans un cône convexe, nous pouvons engendrer les données suivant les lois \sas s et max-stables de façon analogue. Notre méthode ne nécessite pas la discrétisation de la mesure spectrale, elle est donc plus pratique pour simuler les vecteurs aléatoires \sas s et max-stables avec la mesure spectrale absolument continue.  

Nous résumons dans la section \ref{preliminaires} les définitions et les résultats de base concernant les lois \sas s dans un cône convexe. Pour la présentation plus détaillée nous renvoyons le lecteur à l'article \cite{Davydov08}. La section \ref{sect-exemples} présente les lois strictement max-stables et $\alpha$-stables avec des exemples concrets. Nous obtenons l'expression explicite des lois jointes et marginales des vecteurs aléatoires strictement max-stables avec les mesures spectrales données. Les résultats de simulation des vecteurs aléatoires dans les deux familles sont présentés. 
\section{Lois strictement stables dans un cône convexe}\label{preliminaires}
\subsection{Cônes convexes}

\begin{definition}\label{defcone}
Un {\em cône convexe} $\cok$ est un  semigroupe abélien topologique, supposé complet et séparable, avec une opération de multiplication par des scalaires positifs, $(x, a)\rightarrow ax$, continue pour tout $x\in \cok$ et $a>0$ tel que les conditions suivantes sont remplies:

1) $a(x+y)=ax+ay, \ a>0, \ x,y \in \cok$

2) $a(bx)=(ab)x, \ a,b>0, \ x\in \cok$

3) $1x=x, \ x\in\cok$

4) $a{\bf e}={\bf e}, \ a>0, \ \bf{e} \ \mbox{est l'élément neutre de \cok}$.

Le cône $\cok$ s'appelle {\em cône pointé} s'il y a un élément unique {\bf 0} dit {\em l'origine} tel que $ax\rightarrow{\bf 0}$ lorsque $a\downarrow 0$ pour tout $x \in \cok \ \backslash\{{\bf e}\}$.
\end{definition}

Remarquons que la condition suivante qui est typiquement imposée dans la littérature (voir par exemple \cite{Keimel92}) pour un cône n'est pas nécessaire ici
\begin{equation}\label{seconddistri}
(a+b)x=ax+bx, \; a,b>0, \; x\in\cok.
\end{equation}
Par exemple $\Bbb{R}_{+}=[0, \infty)$ avec l'opération d'addition définie par le maximum ne vérifie pas cette condition. Dans un cône pointé, si la distributivité (\ref{seconddistri}) a lieu ou s'il existe l'élément $x\neq{\bf e}$ qui possède l'inverse $(-x)$, alors ${\bf e}={\bf 0}$. Généralement l'élément neutre ne coïncide pas nécessairement avec l'origine.

\begin{definition}
Un élément $z\in\cok$ est {\em $\alpha$-stable} avec $\alpha\neq 0$, si
\[a^{1/\alpha}z+b^{1/\alpha}z=(a+b)^{1/\alpha}z\]
pour tous $a,b>0$.
\end{definition}
Dans la suite on note $\cok(\alpha)$ l'ensemble des éléments $\alpha$-stables de \cok. Il est clair que ${\bf e}, {\bf 0}\in \cok(\alpha)$ pour tout $\alpha\neq 0$. 

\begin{definition}
Le cône \cok \ s'appelle un {\em cône normé} si \cok \ est métrisable par une distance $d$ qui est homogène à l'origine, i.e. $d(ax,{\bf 0})=ad(x,{\bf 0})$ pour tout $a>0$ et $x \in \cok$. La valeur $\|x\|=d(x,{\bf 0})$ s'appelle la {\em norme} de $x$.
\end{definition}

Dans la suite on suppose que \cok \ est un cône normé. Il est clair que $\|x\|=0$ si et seulement si $x={\bf 0}$. De plus si ${\bf e}\neq{\bf 0}$, alors la propriété 4) de la définition \ref{defcone} implique que $\|{\bf e}\|=d({\bf e},{\bf 0})=\infty$. Il est donc essentiel d'admettre que $d$ peut prendre la valeur infinie. Par exemple, si \cok \  est le cône $\overline{\mathbb{R}}_+=[0,\infty]$ où l'opération d'addition est définie par le minimum, noté $(\overline{\mathbb{R}}_+, \wedge)$, alors la distance euclidienne d'un élément non-zéro $x\in\mathbb{R}_+$ à $\infty$ (étant l'élément neutre) est infinie.

L'ensemble
\[S=\{x \; | \; \|x\|=1\}\]
s'appelle la {\em sphère unité}. Le cône pointé admet la {\em décomposition polaire} réalisée par
la bijection $x\leftrightarrow(\|x\|,x/\|x\|)$ entre  
\[\cok '=\cok\backslash \{{\bf 0},{\bf e}\}\] 
et $ (0,\infty)\times S$. 

La distance ou la norme dans \cok \ est dite {\em sous-invariante} si
\begin{equation*}
d(x+h, x)\leq d(h,{\bf 0})=\|h\|, \; x,h\in \cok.
\end{equation*}
Cette contrainte n'est pas triviale, par exemple le cône $(\overline{\mathbb{R}}_+, \wedge)$ avec la distance euclidienne ne vérifie pas cette propriété. Si \cok \ a une norme sous-invariante, alors ${\bf 0}={\bf e}$, et pour tout $\alpha\in (0,1)$, ${\bf e}$ est l'élément unique de norme finie qui appartient à $\cok(\alpha)$.

Les exemples typiques de cônes qui possèdent ces propriétés sont les espaces de Banach ou les cônes convexes dans l'espace de Banach, la famille des sous-ensembles compacts (ou compacts convexes) d'un espace de Banach avec l'addition de Minkowski \cite{Davydov00a} \cite{Gine85} \cite{Gine82}, la famille des ensembles compacts dans $\rd$ avec l'opération d'addition définie par l'union \cite{Molchanov05} et la famille des mesures finies avec l'addition conventionnelle et la multiplication scalaire \cite{Daley03} \cite{Rachev91}. Un autre exemple typique est l'espace $\mathbb{R}_+=[0,\infty)$ avec l'opération d'addition définie par le maximum, i.e. $x+y=x\vee y=\max(x,y)$. D'autres exemples et informations peuvent être trouvés dans \cite{Davydov08}.

\subsection{Processus ponctuels dans un cône}\label{subsect-pp}
Considérons un cône normé \cok \ muni de la $\sigma$-algèbre borélienne $\mathcal{B}(\cok)$. La boule ouverte de rayon $r$ centrée au point {\bf 0} est notée par
\[B_r=\{x\in\cok \; | \; \|x\|<r\}.\]
L'intérieur de son complément est noté par
\[B^r=\{x \; | \; \|x\|>r\}.\] 
Une {\em mesure de comptage}  $m$ est une mesure de forme suivante  
\[m=\delta_{x_1}+\delta_{x_2}+\cdots=\sum_i\delta_{x_i}\]
où $\{x_1,x_2,\ldots\}$ est un ensemble de points dénombrable au maximum. Soit $\mathcal{M}_0$ (respectivement $\mathcal{M}$) la famille des mesures de comptage $m$ sur $\mathcal{B}(\cok)$ telles que $m(B^r)<\infty$ (respectivement $m(B_r)<\infty$) pour tout $r>0$. Pour couvrir ces deux cas avec la même notation, on note $A_r$ au lieu de $B^r$ (respectivement $B_{r}$) dans le cas où l'on considère les mesures dans $\mathcal{M}_0$ (respectivement $\mathcal{M}$). Alors nous avons toujours $m(A_r)<\infty$, c'est-à-dire seulement un nombre fini de points $x_i$ du support de $m$ se situent dans $A_r$. 

Un {\em processus ponctuel $\mu$} est une application mesurable d'un espace de probabilité à $\mathcal{M}_0$ (ou $\mathcal{M}$) muni de la $\sigma$-algèbre engendrée par les ensembles des mesures $m\in\mathcal{M}_0$ (ou $m\in\mathcal{M}$) telle que $m(B)=n$ pour les ensembles boréliens
$B\subset\cok$ et $n\geq 0$. 

Soit $\Lambda$ une mesure sur \cok \ qui est finie sur $A_r$ pour tout $r>0$. Un processus ponctuel $\Pi$ est {\em un processus poissonnien} de mesure d'intensité $\Lambda$, si pour toute famille d'ensembles boréliens disjoints $F_1,\ldots, F_n$, les variables aléatoires $\Pi(F_1),\ldots,\Pi(F_n)$ sont mutuellement indépendantes et distribuées suivant les lois de Poisson d'espérance $\Lambda(F_1),\ldots,\Lambda(F_n)$ respectivement.

\begin{definition}
(Processus ponctuel poissonnien stable $\Pi_{\alpha,\sigma}$) Soit $\Lambda_{\alpha,\sigma}$ une mesure sur $\cok\,'$ définie par $m_\alpha\times\sigma$ où $\sigma$ est une mesure finie sur $\mathcal{B}(S)$ et $m_\alpha$ est une mesure sur $(0,\infty)$ telle que pour $\alpha\neq 0$ on a
\begin{equation}\label{defmalpha}
\left\{
\begin{array}{ll}
m_\alpha((r,\infty))=r^{-\alpha}&\mbox{si} \ \alpha>0,\\
m_\alpha((0,r))=r^{-\alpha}&\mbox{si} \ \alpha<0,
\end{array}\right. r>0.
\end{equation}
Le {\em processus ponctuel poissonnien stable} est un processus poissonnien de mesure d'intensité $\Lambda_{\alpha,\sigma}$, noté $\Pi_{\alpha,\sigma}$.
\end{definition}

La mesure $\sigma$ s'appelle la {\em mesure spectrale}. L'importance du processus $\Pi_{\alpha,\sigma}$ vient du fait que le support de $\Pi_{\alpha,\sigma}$ est union-stable, i.e. $\Pi_{\alpha,\sigma}$ vérifie l'égalité (\ref{stabpi}). 

Pour tout $\alpha\neq 0$, le processus $\Pi_{\alpha,\sigma}$ admet la représentation (\ref{ppintro}) (Th. 3.3 \cite{Davydov08}) avec les éléments aléatoires $\{\epsilon_k\}$  i.i.d. à valeurs dans $S$ et de loi ${\tilde{\sigma}(\cdot)=\sigma(\cdot)/\sigma(S)}$ et $c=\sigma(S)^{1/\alpha}$. 

Soit $\{\xi_k, \ k\geq 1\}$ une suite d'éléments aléatoires i.i.d. à valeurs dans \cok. Un {\em processus ponctuel empirique} est défini par
\begin{equation*}
\beta_n=\sum^{n}_{k=1}\delta_{\xi_k/b_n}, \ n\geq 1,
\end{equation*}
où $\{b_n, \ n\geq 1\}$ est une suite de constantes de normalisation telle que $b_{n}\rightarrow\infty$, si $n\rightarrow\infty$. On va montrer dans le théorème \ref{1.1} que pour avoir de propriétés intéressantes, les constantes $b_n$ doivent avoir la forme suivante
\begin{equation}\label{bn}
b_n=n^{1/\alpha}L(n), \ n\geq 1,
\end{equation}
avec $\alpha\neq 0$ et $L(n)$ une fonction à variation lente à l'infini.

La convergence faible des processus ponctuels $\mu_n\Rightarrow\mu$ a lieu, si $\reie h(\mu_n)\rightarrow \reie h(\mu)$ lorsque $n\rightarrow\infty$ pour toute fonction $h$ sur $\mathcal{M}$ ou $\mathcal{M}_0$ à valeurs dans $\Bbb{R}$ continue dans la topologie vague et bornée. Le processus poissonnien $\Pi_{\alpha,\sigma}$ est la limite faible du processus empirique $\beta_n$, si la loi des $\xi_k$ a la queue à variation régulière \cite{Resnick87}. On donne ici ce résultat sous la forme de \cite{Davydov08}.

\begin{theoreme}(Th. 4.3 \cite{Davydov08})\label{1.1}
Soient $\xi, \xi_1, \xi_2,\ldots$ des éléments aléatoires i.i.d. à valeurs dans $\cok\,'$. Alors $\beta_n\Rightarrow\Pi_{\alpha,\sigma}$ lorsque $n\rightarrow\infty$ pour $\alpha >0$ si et seulement s'il existe une mesure finie $\sigma$ sur $\mathcal{B}(S)$ telle que on a
\begin{equation}\label{regulier2}
\lim_{n\rightarrow\infty}n\reip\left\{\frac{\xi}{\|\xi\|}\in B, \|\xi\|>rb_n\right\}=\sigma(B)r^{-\alpha}, 
\end{equation}
pour tout $r>0$ et $B\in\mathcal{B}(S)$ avec $\sigma(\partial B)=0$, où $b_{n}$ est de forme (\ref{bn}).
\end{theoreme}

\remarque La condition (\ref{regulier2}) est typique pour les théorèmes limites de somme des éléments aléatoires (voir \cite{Resnick87} page 154 et \cite{Araujo80} page 167). Il existe un résultat similaire pour $\alpha<0$ (Cor. 4.4 \cite{Davydov08}). 
\vspace{0.5cm}

La loi d'un élément aléatoire $\xi$ à valeurs dans $\cok\,'$ est dite {\em à queue régulière} si la condition de variation régulière (\ref{regulier2}) est vérifiée. 

\subsection{Série de LePage et domaines d'attraction de loi \sas}\label{subsect-sas}
Un élément aléatoire $X$ dans $\cok$ a la loi {\em strictement $\alpha$-stable} (\sas) s'il vérifie la condition (\ref{stabintro}) de l'introduction avec les opérations d'addition et multiplication définies sur $\cok$. Le théorème ci-dessous donne une famille riche de lois \sas s  par leur décomposition en série.

\begin{theoreme}\label{thlepage}(Th. 3.6 \cite{Davydov08})
Soient $\{\lambda_k, \ k\geq 1\}$ et $\{\epsilon_k, \ k\geq 1\}$ deux suites indépendantes de variables aléatoires i.i.d.. Les variables aléatoires $\lambda_k$ ont la loi exponentielle standard, les éléments aléatoires $\epsilon_k$ sont de même loi $\tilde{\sigma}(\cdot)$ qui est une mesure finie normalisée sur $S$. Notons $\Gamma_k=\lambda_1+\cdots+\lambda_k, \ k\geq1$. Si la valeur principale de l'intégrale $\int x\Pi_{\alpha,\tilde{\sigma}}(dx)$ est finie, i.e. $\int_{A_{r}}x\Pi_{\alpha,\tilde{\sigma}}(dx)<\infty$, avec probabilité $1$, alors pour tout $z\in\cok(\alpha)$ et $c\geq0$, la série
\begin{equation}\label{serlepage}
\xi_{\alpha,\sigma}=z+c\sum_{k=1}^\infty\Gamma_k^{-1/\alpha}\epsilon_k
\end{equation}
converge presque sûrement (p.s.) et $\xi_{\alpha,\sigma}$ admet une loi \sas \ sur \cok. Le paramètre $\alpha$ s'appelle {\em l'indice caractéristique} et la mesure $\sigma(\cdot)=c^{\alpha}\tilde{\sigma}(\cdot)$ s'appelle {\em la mesure spectrale}.

Si la norme de \cok \ est sous-invariante, alors la somme infinie (\ref{serlepage}) converge absolument p.s. pour tout $\alpha\in (0,1)$.
\end{theoreme}

\begin{definition}
Un élément aléatoire $X$ à valeurs dans \cok \ appartient {\em au domaine d'attraction} d'un élément aléatoire \sas, noté $\xi_{\alpha,\sigma}$, si pour $\{X_n, \ n\geq 1\}$, une suite de copies indépendantes de $X$, on a
\begin{equation}\label{dom}
b_n^{-1}(X_1+\cdots +X_n)\Rightarrow\xi_{\alpha,\sigma},
\end{equation} 
où $\{b_n,\ n\geq 1\}$ est une suite de constantes de normalisation positives et $\Rightarrow$ désigne la convergence faible des éléments aléatoires à valeurs dans \cok.
\end{definition}

Le fait que l'élément aléatoire $X$ appartienne au domaine d'attraction de $\xi_{\alpha,\sigma}$ sera noté par l'écriture ``$X\in \mbox{Dom}(\alpha,\sigma)$''. Le théorème suivant montre qu'avec une condition supplémentaire sur \cok \  la régularité de la queue de la loi de $X$ implique que $X\in \mbox{Dom}(\alpha,\sigma)$ pour $0<\alpha<1$.

\begin{theoreme}\label{4.7}(Th. 4.7 \cite{Davydov08})
Supposons que $\cok$ \ a la norme sous-invariante. Si la loi d'un élément aléatoire $X\in \cok\,'$ vérifie la condition (\ref{regulier2}) avec $\alpha\in(0,1)$, alors $X$ appartient au domaine d'attraction d'un élément aléatoire \sas \, qui admet la représentation du type de LePage (\ref{lepageintro}) avec $\Gamma_k$, $\epsilon_k$ et $c$ qui sont définis comme dans le théorème \ref{thlepage}. 
\end{theoreme}

\remarque La condition (\ref{regulier2}) est nécessaire et suffisante pour qu'un vecteur aléatoire à valeurs dans $\rd$ appartienne au domaine d'attraction d'une loi \sas \, avec $0<\alpha<2$ \cite{Araujo80}.
\vspace{0.5cm}

L'existence de la mesure spectrale d'une loi \sas \, dans un espace de Banach séparable est connue. L'existence de la représentation de LePage d'un élément aléatoire \sas \, dans un semigroupe plus général est démontrée dans \cite{Davydov08}. Dans la section suivante nous présentons quelques exemples de lois max-stables et $\alpha$-stables. 
\section{Quelques exemples et simulation}\label{sect-exemples}

\subsection{Lois strictement max-stables}\label{subsect-paramax}

On considère les lois stables dans le cône $(\mathbb{R}^d_{+}, \vee)$. Les vecteurs dans $\mathbb{R}^d_{+}=[0,\infty)^d$ sont notés par $\reixs=(x_1,\ldots, x_d)$. Les relations et opérations sont définies par composante, c'est-à-dire pour $\reixs , \reiys\in\mathbb{R}^d_{+} $
\begin{eqnarray*}
\reixs< \reiys & \mbox{signifie} &  x_i< y_i,  i=1,\ldots,d,\\
\reixs\leqslant \reiys  & \mbox{signifie} & x_i\leqslant y_i, i=1,\ldots,d,\\
\reixs+\reiys &: =&\reixs\vee \reiys  =(x_1\vee y_1,\ldots, x_d\vee y_d),\\
a\reixs &: = & (ax_1,\ldots,ax_d), \;\; a>0.
\end{eqnarray*}
Les rectangles sont notés par 
\[(\reia, \reib)=\{(x_1,\ldots, x_d) \;|\; a_1< x_1< b_1,\ldots, a_d< x_d < b_d\},\]
\[[\reia, \reib]=\{(x_1,\ldots, x_d) \;|\; a_1\leqslant x_1\leqslant b_1,\ldots, a_d\leqslant x_d \leqslant b_d\}.\]
Prenons la $L_{\infty}$-norme, i.e. $\|\reixs\|=\max (x_1,\ldots, x_d)$. La sphère unité dans cette norme $\{\reixs\; |\; \|\reixs\|=1\}$ est notée par $S_{\vee}^{d-1}$. 

Un vecteur aléatoire $X=(X^{(1)},\ldots,X^{(d)})$ a la loi {\em strictement max-stable} si pour tous $a,b>0$
\[\left(a^{1/\alpha}X_{1}^{(1)}\bigvee b^{1/\alpha}X_{2}^{(1)},\ldots,a^{1/\alpha}X_{1}^{(d)}\bigvee b^{1/\alpha}X_{2}^{(d)}\right)\egenloi (a+b)^{1/\alpha}(X^{(1)},\ldots,X^{(d)})\] 
où $X_{1}$, $X_{2}$ sont les copies indépendantes de $X$. Soit $\epsilon_1,\epsilon_2,\ldots$ une suite de vecteurs aléatoires i.i.d. à valeurs dans $S^{d-1}_\vee$ suivant la loi $\sigma$. Soit $\lambda_{1},\lambda_{2},\ldots$ une suite de variables aléatoires i.i.d. de loi exponentielle standard et $\Gamma_i=\sum_{j=1}^{i}\lambda_{j}$. Les suites $\{\epsilon_i\}$ et $\{\Gamma_i\}$ sont indépendantes. La série suivante converge absolument p.s. pour tout $\alpha>0$
\begin{equation}\label{lepagemax}
\bigvee_{i=1}^\infty \Gamma_i^{-1/\alpha}\epsilon_i.
\end{equation}
Le résultat direct du théorème \ref{thlepage} indique que cette série produit un vecteur aléatoire strictement max-stable avec l'indice $\alpha>0$ et la mesure spectrale $\sigma$. Dans la suite on utilise l'écriture ``$\mathcal{MS}_{d}(\alpha,\sigma)$'' pour indiquer la loi strictement max-stable $d$-dimensionnelle d'indice $\alpha$ et de mesure spectrale $\sigma$.

Soit $\xi_{\alpha,\sigma}$ un vecteur aléatoire de loi $\mathcal{MS}_{d}(\alpha,\sigma)$. Puisque $\xi_{\alpha,\sigma}$ et la série (\ref{lepagemax}) ont la même loi, le maximum partiel suivant donne une approximation de la loi de $\xi_{\alpha,\sigma}$
\[\hat{\xi}_{k}=\bigvee_{i=1}^k \Gamma_i^{-1/\alpha}\epsilon_i.\]
Notons
\[\epsilon_{i}=(\epsilon_{i}^{(1)},\ldots,\epsilon_{i}^{(d)})\;\;\;\mbox{et}\;\;\;\tau_{j}=\min\{i ~|~ \epsilon_{i}^{(j)}=1\}, \ j=1,\ldots,d,\, i=1,2,\ldots.\]
Supposons que
\[\reip\{\epsilon_1^{(j)}=1\}=p_{j}> 0, \; j=1,\ldots,d.\]
Le vecteur aléatoire $\hat{\xi}_{k}$ converge vers $\xi_{\alpha,\sigma}$ en probabilité lorsque $k$ tend vers l'infini. En effet
\[\reip\{\hat{\xi}_{k}\neq\xi_{\alpha,\sigma}\}=\reip\{\exists j \; \mbox{tel que} \; \tau_j>k\}\leq \sum_{j=1}^{d}\reip\{\tau_j>k\}=\sum_{j=1}^{d}(1-p_{j})^{k}.\]
La dernière égalité vient du fait que la loi de $\tau_j$ est géométrique
\[\reip\{\tau_j=n\}=\reip\{\epsilon_1^{(j)}<1,\ldots,\epsilon_{n-1}^{(j)}<1,\epsilon_n^{(j)}=1\}=(1-p_{j})^{n-1}p_{j}.\]

Notons le processus ponctuel
\begin{equation}\label{maxppp}
\Pi_{\alpha,\sigma}(\cdot)=\sum_{i=1}^\infty\delta_{\Gamma_{i}^{-1/\alpha}\epsilon_{i}}(\cdot),
\end{equation}
où $\Gamma_{i}$ et $\epsilon_{i}$ sont définis comme dans (\ref{lepagemax}). D'après la représentation du processus ponctuel poissonnien stable (\ref{ppintro}) et la définition (\ref{defmalpha}), la mesure d'intensité de $\Pi_{\alpha,\sigma}$ est $m_\alpha\times\sigma$. Remarquons que pour tout $r>0$, puisque $m_\alpha((r,\infty))=r^{-\alpha}<\infty$, on a
\[\reip\{\Pi_{\alpha,\sigma}(\{\reixs\,|\,\|\reixs\|>r\})<\infty\}=1.\]
Il y a seulement un nombre fini de points du support de $\Pi_{\alpha,\sigma}$ qui sont grands. Le fait que $\xi_{\alpha,\sigma}$ et la série (\ref{lepagemax}) ont la même loi implique que la loi de $\xi_{\alpha,\sigma}$ ne dépend que des points les plus grands du processus $\Pi_{\alpha,\sigma}$. La loi de $\xi_{\alpha,\sigma}$ est facile à calculer, pour $\reixs>0$
\begin{equation}\label{jointmax}
\reip\{\xi_{\alpha,\sigma}\in[0,\reixs]\}=\reip\{\Pi_{\alpha,\sigma}([0,\reixs]^{\complement})=0\}=\exp(-m_{\alpha}\times\sigma([0,\reixs]^{\complement})).
\end{equation}

La représentation de LePage (\ref{lepagemax}) peut être écrite sous la forme suivante
\[\xi_{\alpha,\sigma}\egenloi\left(\bigvee_{k=1}^\infty\Gamma_{k}^{-1/\alpha}\epsilon_{k}^{(1)},\ldots,\bigvee_{k=1}^\infty\Gamma_{k}^{-1/\alpha}\epsilon_{k}^{(d)}\right)=(\Gamma^{-1/\alpha}_{\tau_{1}},\ldots,\Gamma^{-1/\alpha}_{\tau_{d}}).\]
La loi marginale de $\xi_{\alpha,\sigma}$ peut être obtenue en utilisant le processus ponctuel  
\begin{equation}\label{maxpppmarg}
\Pi_{\alpha,\sigma^{(j)}}=\sum_{i=1}^\infty\delta_{\Gamma_{i}^{-1/\alpha}\epsilon_{i}^{(j)}}, \ j=1,\ldots,d.
\end{equation}
Si la loi de la composante $\epsilon_{1}^{(j)}$, notée $\sigma^{(j)}$, est connue, alors
\begin{equation}\label{margmax}
\reip\{\xi_{\alpha,\sigma}^{(j)}\leq x_j\}=\reip\{\Pi_{\alpha,\sigma^{(j)}}((x_j,\infty))=0\}=\exp(-m_{\alpha}\times\sigma ^{(j)} ((x_j,\infty))).
\end{equation}

\vspace{0.5cm}

Dans les trois exemples ci dessous nous utilisons le lien entre le processus ponctuel poissonnien et la loi strictement max-stable pour déduire la fonction de répartition de la loi avec la mesure spectrale donnée. L'exemple \ref{ex1} présente la mesure spectrale discrète et concentrée sur les points d'intersection des axes et de la sphère unité. Dans ce cas le vecteur aléatoire strictement max-stable a ses composantes indépendantes. L'exemple \ref{ex2} présente la mesure spectrale uniforme. L'expression explicite de la loi jointe est obtenue pour $d=2$. L'exemple \ref{ex3} présente la mesure spectrale  dans $\Bbb{R}^2$ concentrée sur deux points qui ne sont pas nécessairement sur les axes. L'exemple \ref{ex1} pour $d=2$ est un cas particulier de l'exemple \ref{ex3}. Des vecteurs aléatoires strictement max-stables à valeurs dans $\Bbb{R}^2$ avec les mesures spectrales présentées sont simulés en utilisant le maximum partiel de la représentation de LePage $\hat{\xi}_{k}$. 

\vspace{0.5cm}


\exemple{ex1} Considérons le vecteur aléatoire $\xi_{\alpha,\sigma}$ de loi $\mathcal{MS}_{d}(\alpha,\sigma)$ avec
\[\sigma(\cdot)=\sum\limits_{i=1}^{d} p_{i}\delta_{e_{i}}(\cdot),\] 
où $e_{i}=(e_{i}^{(1)},\ldots,e_{i}^{(d)})$, $e_{i}^{(j)}=\left\{\begin{array}{ll}1,&i=j\\0,&\mbox{sinon}\end{array}\right.$ et $\sum\limits_{i=1}^{d}p_{i}=1, \ p_{i}\geq 0, \ i=1,\ldots,d$.
C'est-à-dire la mesure spectrale $\sigma$ concentrée sur les points où les axes et la sphère unité se croisent. 

Pour chaque $i$, $\epsilon_{1}^{(i)}$ est la variable aléatoire de loi de Bernoulli de paramètre $p_i$. En appliquant (\ref{margmax}), on a la loi marginale de $\xi_{\alpha,\sigma}$
\[\reip\{\xi_{\alpha,\sigma}^{(i)}\leq x_i\}=\exp(-p_{i}x_i^{-\alpha}).\]

En appliquant (\ref{jointmax}), on a la fonction de répartition de $\xi_{\alpha,\sigma}$, pour $\reixs>0$
\begin{equation*}
\reip\{\xi_{\alpha,\sigma}\in[0,\reixs]\}=\exp\left(-\sum_{i=1}^{d}p_{i}x_i^{-\alpha}\right).
\end{equation*} 

\begin{figure}[!htbp]
\begin{center} 
\begin{tabular}{cc}
a. \includegraphics[width=50mm, height=50mm]{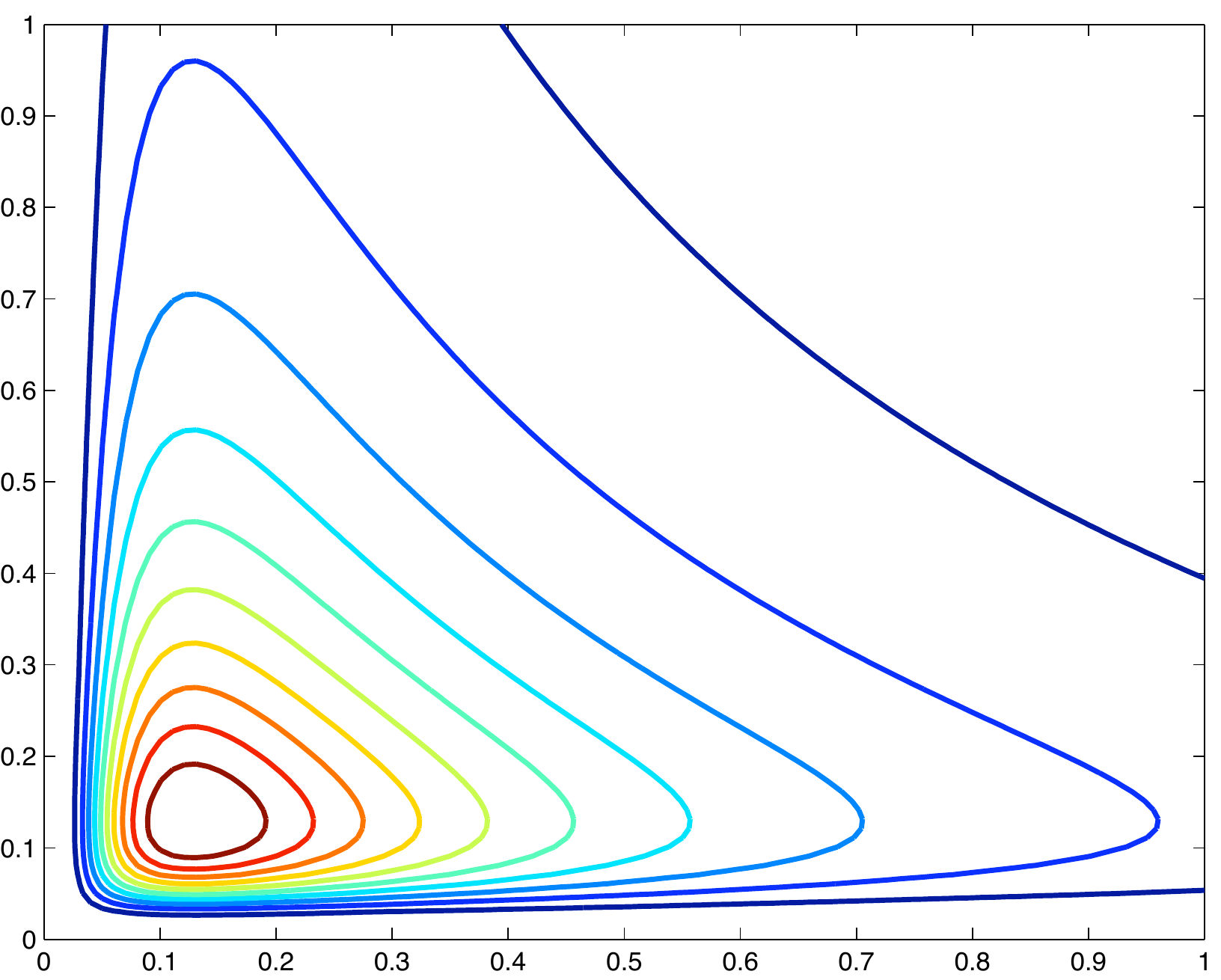} & 
 \includegraphics[width=50mm, height=50mm]{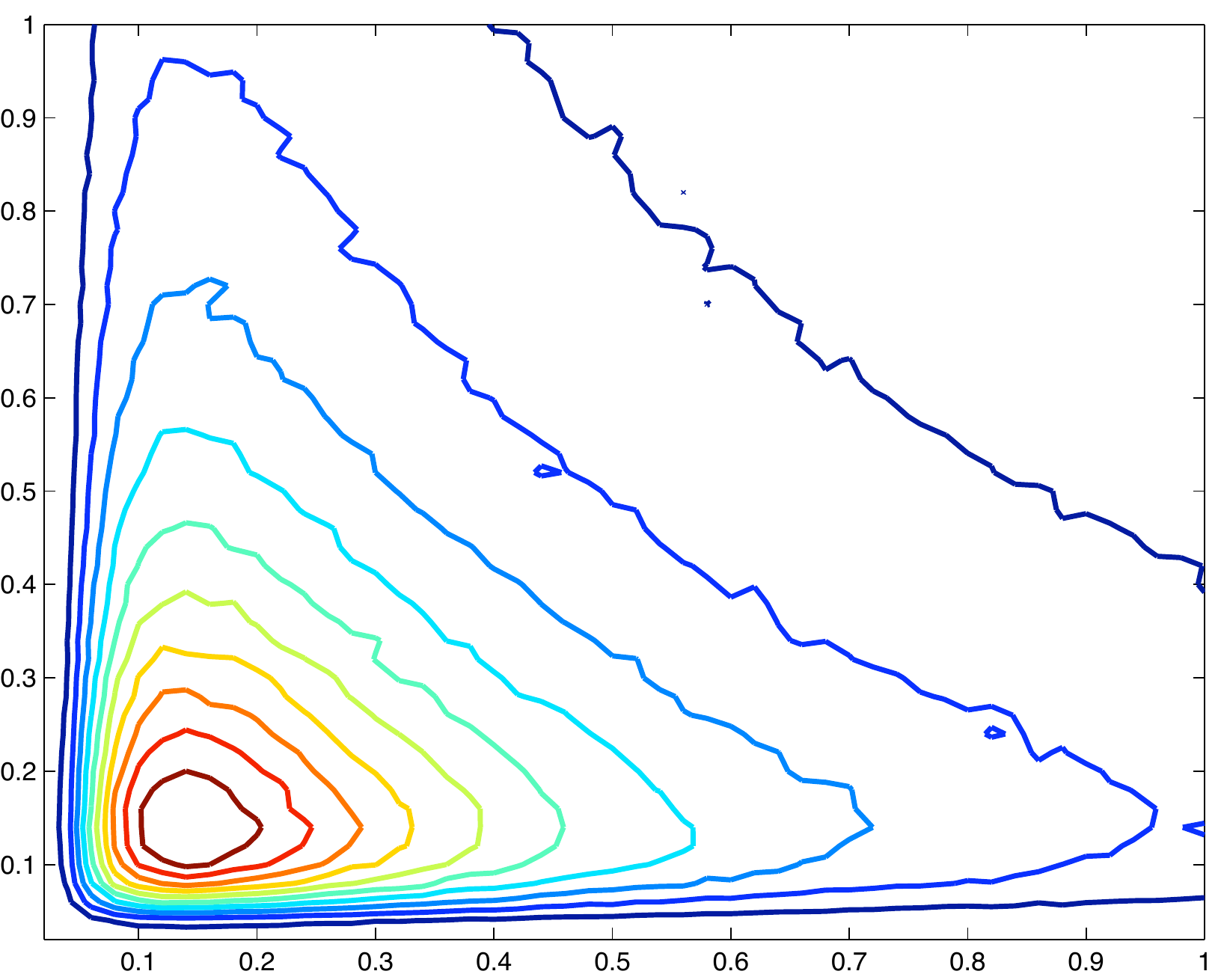} \\
b. \includegraphics[width=50mm, height=50mm]{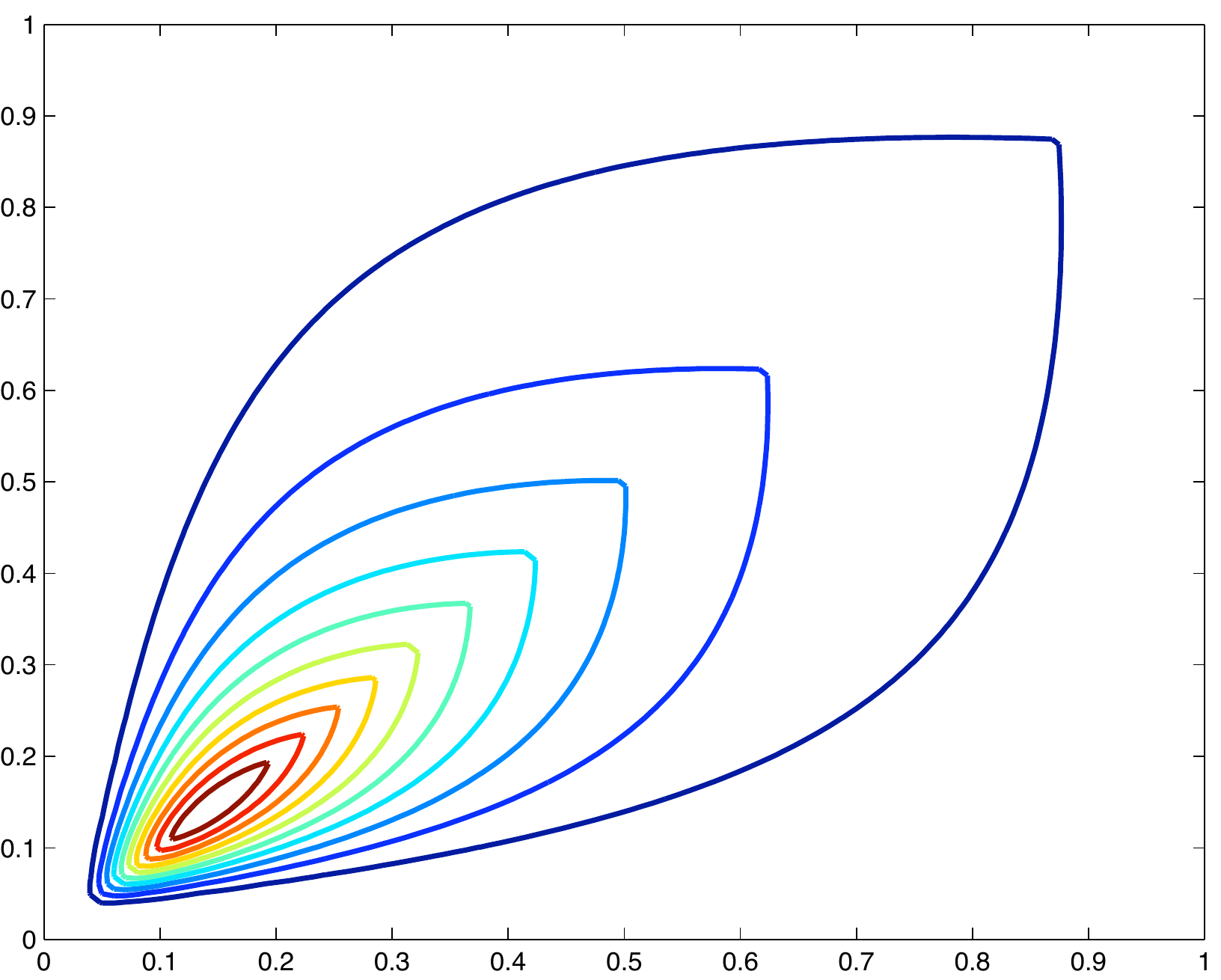} &
\includegraphics[width=50mm, height=50mm]{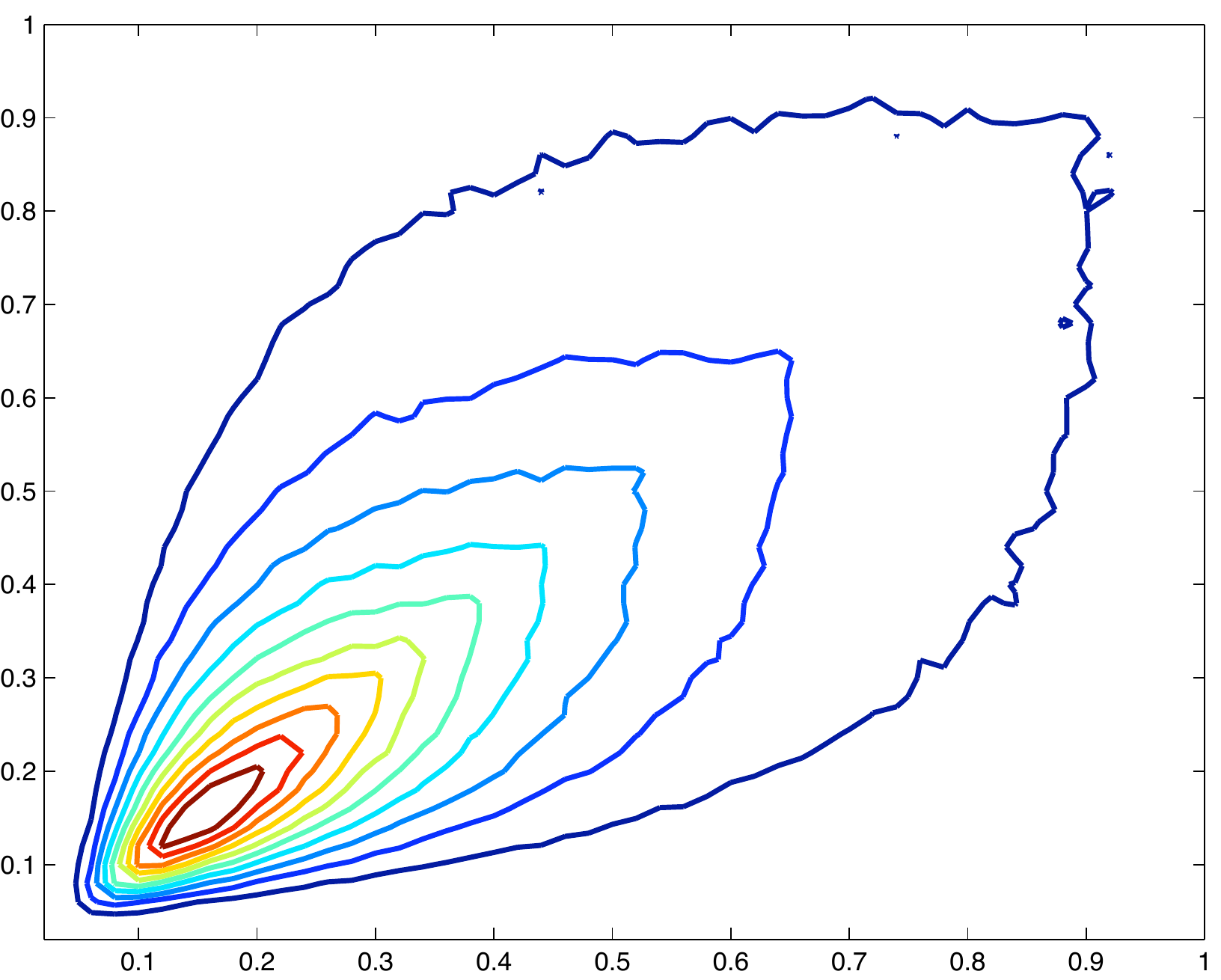} 
\end{tabular}
\end{center}
\caption[Simulations des vecteurs aléatoires de loi $\mathcal{MS}_{2}(0.75,\sigma)$.]{Courbes de niveau des densités réelles (gauche) et empiriques (droite) de la loi $\mathcal{MS}_{2}(0.75,\sigma)$, où $\sigma$ est discrète définie dans l'exemple \ref{ex1} (a.), uniforme définie dans l'exemple \ref{ex2} (b.). Le nombre de données simulées est $10^7$.}
\label{contourmax12}
\end{figure}

En conclusion, le vecteur aléatoire strictement max-stable $\xi_{\alpha,\sigma}$ a ses composantes indépendantes si la mesure spectrale est concentrée sur l'intersection des axes et de la sphère unité. C'est un résultat connu qu'on peut trouver par exemple dans \cite{Resnick87} (Cor. 5.25). La figure \ref{contourmax12} (a.) présente les courbes de niveau de la densité de loi strictement max-stable bivariée d'indice $\alpha=0.75$ et  mesure spectrale $\sigma(\cdot)=0.5\delta_{\{(1,0)\}}(\cdot)+0.5\delta_{\{(0,1)\}}(\cdot)$. 
\vspace{0.5cm}


\exemple{ex2} Soit $\sigma$ la mesure uniforme sur $S_{\vee}^{d-1}$. Soit $\epsilon$ le vecteur aléatoire unitaire de loi $\sigma$. Avant de considérer la loi $\mathcal{MS}_{d}(\alpha,\sigma)$, nous montrons que la loi de composante $\epsilon^{(i)}$, $i=1,\ldots,d$, est
\begin{equation}\label{defsigmaunif}
\sigma^{(i)}(\cdot)=\frac{1}{d}\one_{\{1\}}(\cdot)+\frac{d-1}{d}\mu(\cdot)\one_{\mathcal{B}([0,1))}(\cdot),
\end{equation}
où $\mu$ est la mesure de Lebesgue sur $[0,1)$. Puisque $\|\epsilon\|=1$, on a
\begin{equation}\label{unif1}
\reip\{\exists i\in\{1,\ldots,d\} \; \mbox{tel que} \; \epsilon^{(i)}=1\}=1.
\end{equation}
L'absolue continuité de la loi de $\epsilon$ implique que
\begin{equation}\label{unif2}
\reip\{\epsilon^{(i)}=\epsilon^{(j)}=1, i\neq j\}=0.
\end{equation}
En considérant (\ref{unif1}), (\ref{unif2}) et l'uniformité de la loi de $\epsilon$, on a
\begin{equation}\label{unif3}
\reip\{\epsilon^{(i)}=1\}=\frac{1}{d}, \; i=1,\ldots,d.
\end{equation}
Ainsi on a
\begin{equation}\label{unif4}
\reip\{\epsilon^{(i)}\in [0,1)\}=\frac{d-1}{d}, \; i=1,\ldots,d.
\end{equation}
Les égalités (\ref{unif3}), (\ref{unif4}) et l'absolue continuité de la loi de $\epsilon$ impliquent (\ref{defsigmaunif}).

Soit $\xi_{\alpha,\sigma}$ le vecteur aléatoire de loi $\mathcal{MS}_{d}(\alpha,\sigma)$ avec $\sigma$ uniforme. En appliquant (\ref{margmax}), on a la loi marginale de $\xi_{\alpha,\sigma}$
\begin{equation*}
\reip\{\xi_{\alpha,\sigma}^{(i)}\leq x_i\}=\exp\left(-\frac{\alpha+d}{(\alpha+1)d}x_i^{-\alpha}\right).
\end{equation*}
En effet en notant $\eta$ la variable aléatoire de loi $m_{\alpha}$, on a
\begin{eqnarray*}
m_{\alpha}\times\sigma^{(i)}((x_i,\infty))&=&\reip\{\eta\epsilon_{k}^{(i)}>x_i\}\\
&=&\reip\{\eta>x_i, \epsilon_{k}^{(i)}=1\}+\reip\{\eta\epsilon_{k}^{(i)}>x_i, \epsilon_{k}^{(i)}\in [0,1)\}\\
&=&\frac{1}{d}(x_i)^{-\alpha}+\int_{x_i}^{\infty}\int_{x}^{\infty}\alpha t^{-\alpha-2}\frac{d-1}{d} dt dx\\
&=&\frac{\alpha+d}{(\alpha+1)d}x_i^{-\alpha}.
\end{eqnarray*}
De manière analogue on peut déduire de (\ref{jointmax}) la fonction de répartition de $\xi_{\alpha,\sigma}$. Il n'y a pas d'expression générale pour tout $d$. Dans le cas où $d=2$ la loi jointe de $\xi_{\alpha,\sigma}$ est
\[F_{\xi_{\alpha,\sigma}}(x_1, x_2)=\exp\left(-\frac{1}{2(\alpha+1)}(\alpha x_{(1)}x_{(2)}^{-\alpha-1}+(\alpha+2)x_{(1)}^{-\alpha})\right),\]
où $x_{(1)}=x_1\bigwedge x_2$ et $x_{(2)}=x_1\bigvee x_2$. La figure \ref{contourmax12} (b.) présente les courbes de niveau de la densité de loi strictement max-stable bivariée d'indice $\alpha=0.75$ et de mesure spectrale uniforme. 

On discute finalement la régularité de la queue de cette loi. En prenant l'ensemble $B=[(1,0),(1,t)]$ ou $B=[(0,1),(t,1)], 0\leqslant t\leqslant 1$, on a pour tout $r>0$
\begin{eqnarray*}
I_1&=&\reip\{\Pi_{\alpha,\sigma}((r,\infty)\times B)\neq 0, \Pi_{\alpha,\sigma}((tr,\infty)\times B^\complement)=0\}\\
&\leq&\reip\{\xi_{\alpha,\sigma}\in (r,\infty)\times B\}\\
&\leq&\reip\{\Pi_{\alpha,\sigma}((r,\infty)\times B)\neq 0\}\\
&=&I_2,
\end{eqnarray*}
où \[I_1=(1-\exp(-\sigma(B)r^{-\alpha}))\exp(-\sigma(B^\complement)(tr)^{-\alpha}) \;\;\; \mbox{et} \;\;\; I_2=1-\exp(-\sigma(B)r^{-\alpha}).\]
On en déduit
\begin{equation*}
\reip\{\xi_{\alpha,\sigma}\in B\times (r,\infty)\}\sim \sigma(B)r^{-\alpha} \;\;\;  \mbox{lorsque} \;\;\;  r\rightarrow\infty.
\end{equation*}
Donc la loi de $\xi_{\alpha,\sigma}$ vérifie la condition de variation régulière (\ref{regulier2}).
\vspace{0.5cm}

\begin{figure}[!htbp]
\begin{center} 
\begin{tabular}{cc}
a. \includegraphics[width=50mm, height=50mm]{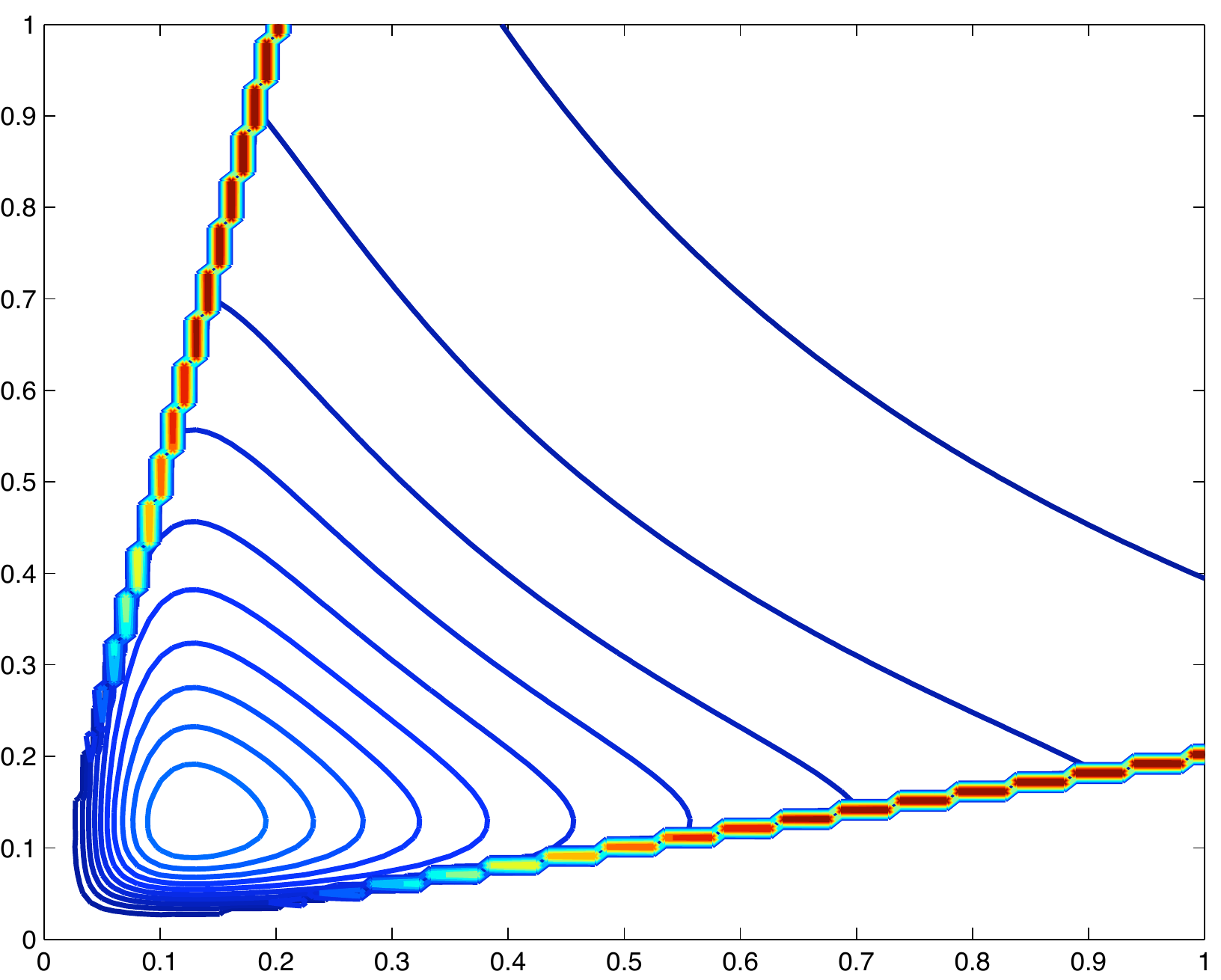} & 
 \includegraphics[width=50mm, height=50mm]{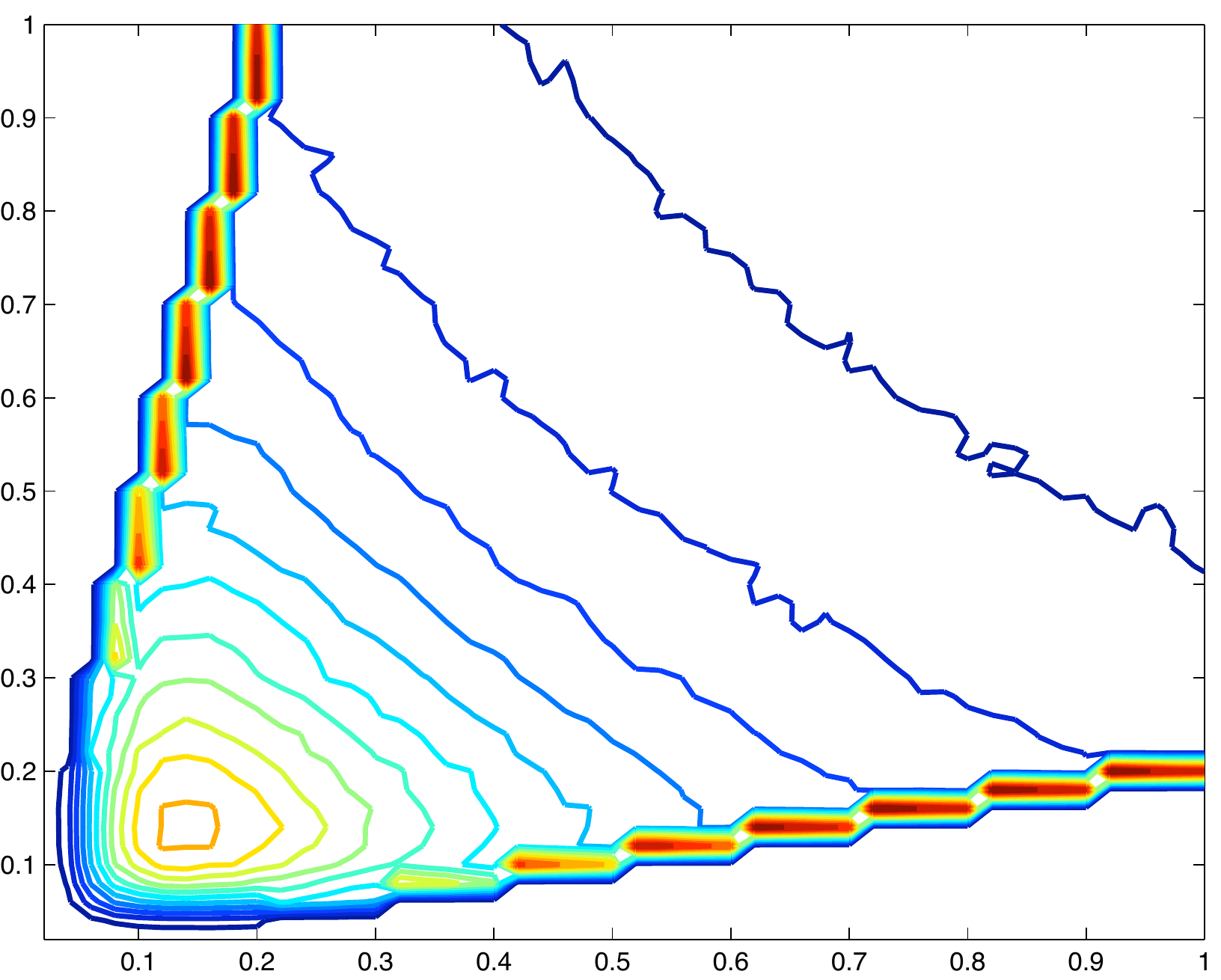} \\
b. \includegraphics[width=50mm, height=50mm]{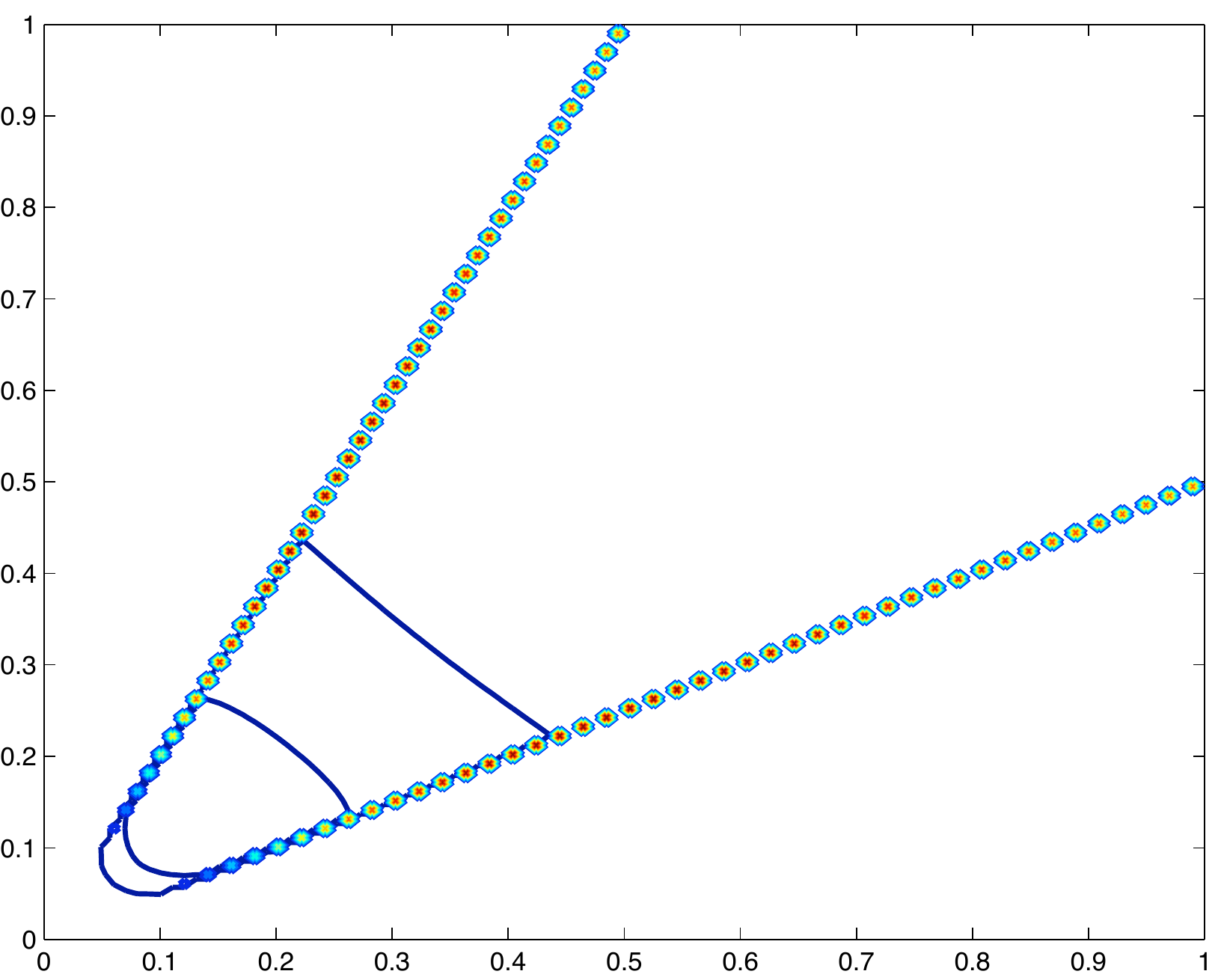} & 
 \includegraphics[width=50mm, height=50mm]{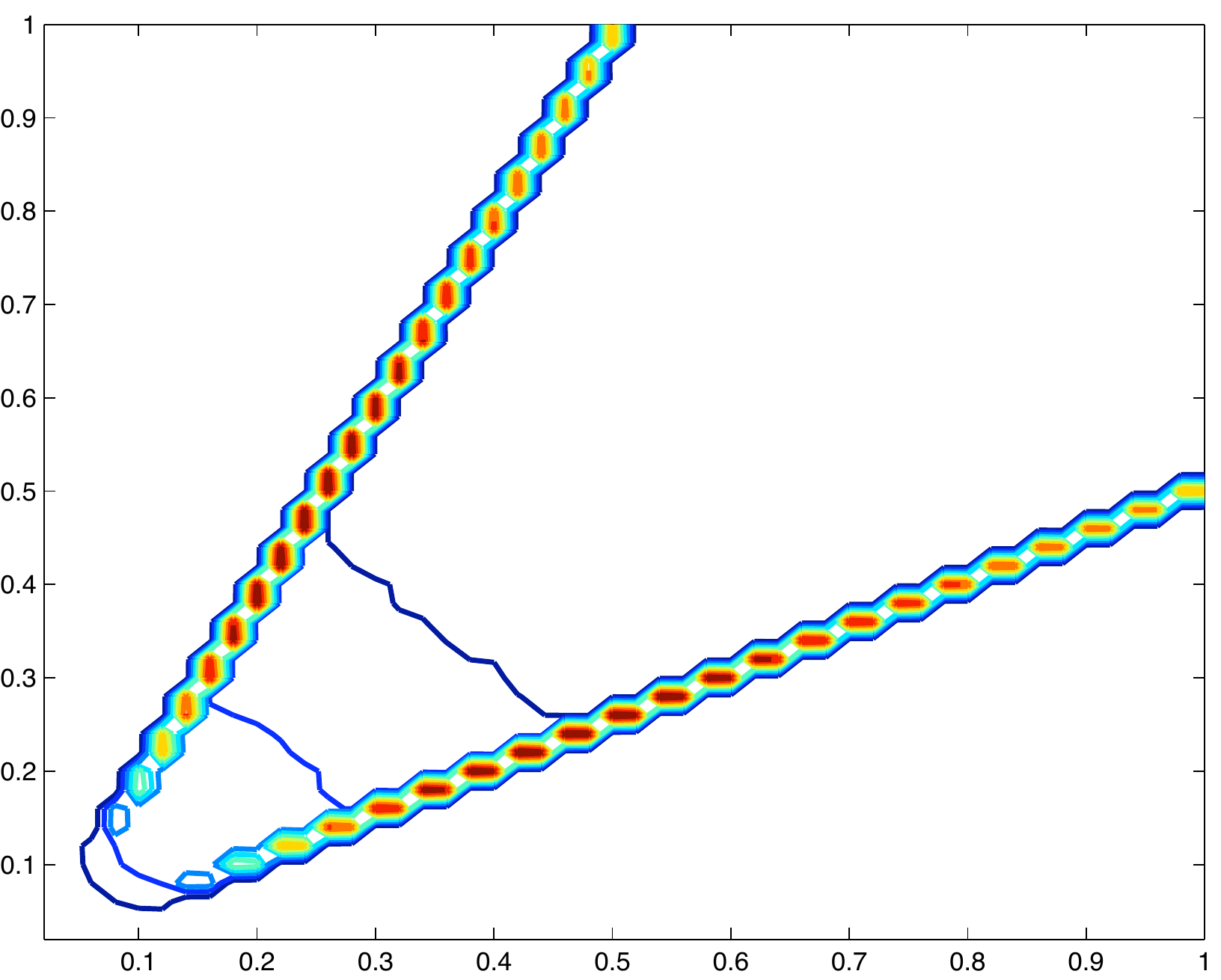}\\
c. \includegraphics[width=50mm, height=50mm]{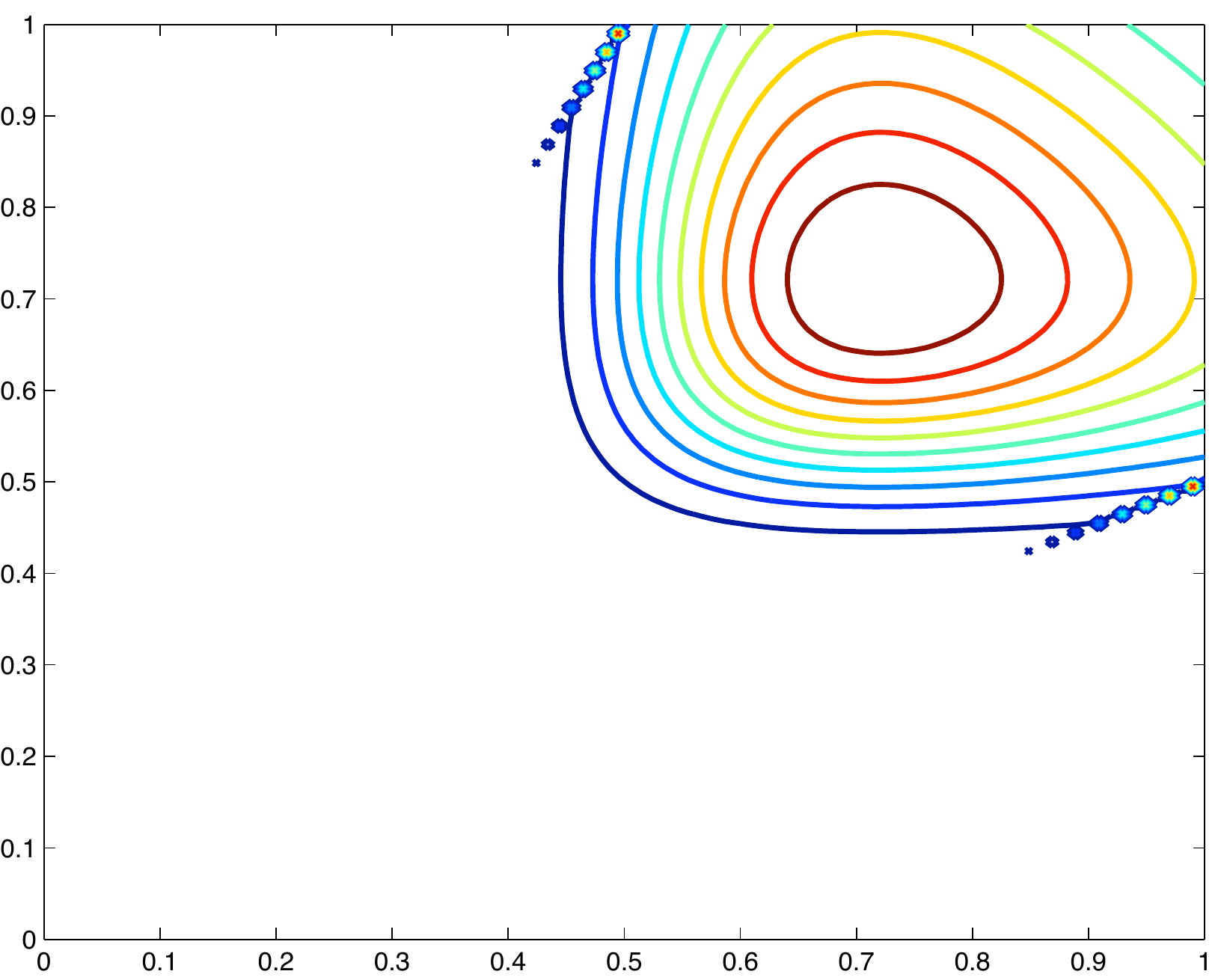} &
 \includegraphics[width=50mm, height=50mm]{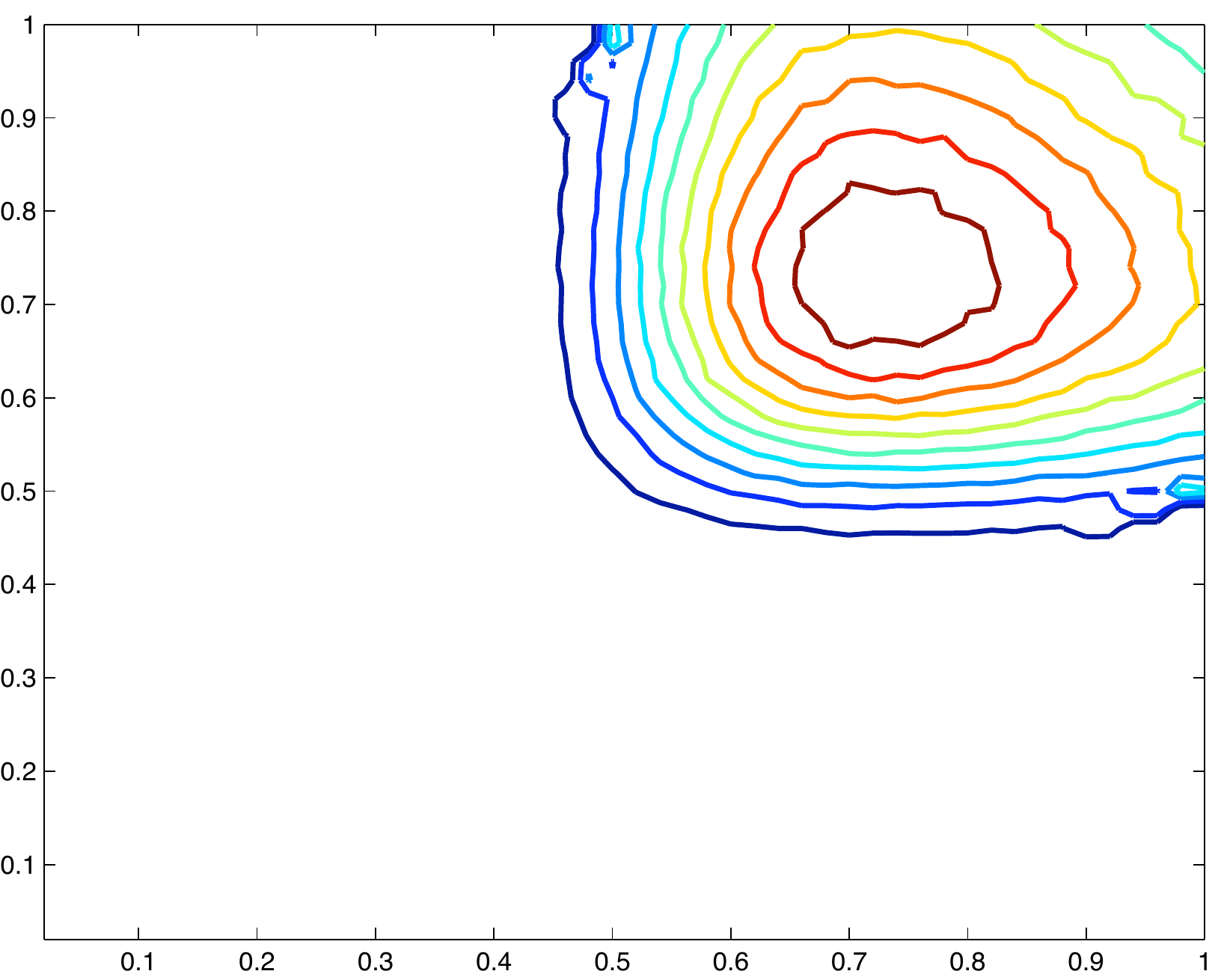}
\end{tabular}
\end{center}
\caption[Simulations des lois $\mathcal{MS}_{2}(\alpha,\sigma)$.]{Courbes de niveau des densités réelles (gauche) et empiriques (droite) de la loi $\mathcal{MS}_{2}(\alpha,\sigma)$ où $\sigma$ est discrète définie dans l'exemple \ref{ex3} avec $p=q=0.5$, $\alpha=0.75$ et $a=b=0.2$ (a.), $\alpha=0.75$ et $a=b=0.5$ (b.), $\alpha=3$ et $a=b=0.5$. Le nombre de données simulées est $10^7$.}
\label{contourmax1}
\end{figure}

\exemple{ex3} Considérons le vecteur aléatoire $\xi_{\alpha,\sigma}$ de loi $\mathcal{MS}_{2}(\alpha,\sigma)$ avec $\sigma$ concentrée sur deux points
\begin{equation}\label{defsigma2p}
\sigma(\cdot)=p\delta_{e_{1}}(\cdot)+q\delta_{e_{2}}(\cdot),
\end{equation}
où $e_{1}=(1,a)$, $e_{2}=(b,1)$, $0\leq a,b \leq 1$ et $p+q=1$. Il est facile de voir que $\xi_{\alpha,\sigma}$ est à valeurs dans la région entre deux demi-droites sortant du point d'origine vers les points $e_{1}$ et $e_{2}$. On en déduit la fonction de répartition 
\[\reip\{\xi\in[0,\reixs]\}=\left\{ \begin{array}{ll}
\exp(-px_1^{-\alpha}-qx_2^{-\alpha}),&\; a<x_2/x_1<1/b,\\
\exp(-(pa^{\alpha}+q)x_2^{-\alpha}),&\; x_2/x_1<a,\\
\exp(-(p+qb^\alpha)x_1^{-\alpha}),&\; x_2/x_1>1/b.\end{array}\right.\]
En notant le moment d'arrêt $\tau=\min\{i\,|\,\epsilon_{i}=e_{2}\}$, considérons la probabilité que $\xi_{\alpha,\sigma}$ est à valeurs dans la demi-droite $\{e_{1}\times (0,\infty)\}$
\begin{eqnarray*}
\reip\{\xi_{\alpha,\sigma}\in e_{1}\times (0,r)\}&=&\sum_{i=2}^\infty\reip\{\epsilon_{1}=e_{1}, \tau=i,\Gamma_{i}^{-1/\alpha}<a\Gamma_{1}^{-1/\alpha}, \Gamma_{1}^{-1/\alpha}<r\}\\
&=&\sum_{i=2}^{\infty}p^{i-1}q\reip\left\{\frac{\lambda_{2}+\cdots+\lambda_{i}}{\lambda_{1}}>a^{-\alpha}-1, \lambda_{1}>r^{-\alpha}\right\}\\
&=&\left(p-\frac{pq(1-a^\alpha)}{q+pa^\alpha}\right)\exp(-(p+qa^{-\alpha})r^{-\alpha}).
\end{eqnarray*}
De manière analogue on a
\[\reip\{\xi_{\alpha,\sigma}\in e_{2}\times (0,r)\}=\left(q-\frac{pq(1-b^\alpha)}{p+qb^\alpha}\right)\exp(-(q+pb^{-\alpha})r^{-\alpha}).\]
Si $a=b=1$, alors $p=q=1$ et la loi de $\xi_{\alpha,\sigma}$ est définie sur la demi droite qui divise le premier quadrant en deux moitiés égales. Si $a=b=0$ alors $\xi_{\alpha,\sigma}$ a ses valeurs strictement à l'intérieur du premier quadrant et ses composantes sont indépendantes. Cela coïncide avec l'exemple \ref{ex1}. Plus $a$ et $b$ sont proches de $1$, plus il est possible que $\xi_{\alpha,\sigma}$ est à valeurs sur les demi-droites. La figure \ref{contourmax1} présente la dépendance des paramètres $a$, $b$ et $\alpha$ de la densité de $\mathcal{MS}_{2}(\alpha,\sigma)$ avec $\sigma$ définie par (\ref{defsigma2p}).

Considérons la régularité de la queue de cette loi. Définissons deux vecteurs aléatoires 
\[\epsilon^{(1)}\egenloi p\delta_{e_{1}}, \; \; \epsilon^{(2)}\egenloi q\delta_{e_{2}}.\]
Notons $\epsilon_{1}^{(i)}, \epsilon_{2}^{(i)}, \ldots$ les copies indépendantes de $\epsilon^{(i)}$, $i=1,2$. D'après le résultat de l'amincissement du processus (voir par exemple Prop. 4.4.1 \cite{Resnick92}), les processus $\Pi_1$ et $\Pi_2$ définis par
\[\Pi_{i}=\sum_{j=1}^\infty\delta_{\Gamma_{j}^{-1/\alpha}\epsilon_{j}^{(i)}}, \;\; i=1,2,\]
sont indépendants et poissonniens. Les mesures d'intensité sont respectivement $m_\alpha\times p\delta_{e_{1}}$ et $m_\alpha\times q\delta_{e_{2}}$. Prenons l'ensemble 
$B=[(1,s), (1,t)], \;\; 0\leq s<a<t\leq 1.$ 
Il est facile de voir que pour tout $ r>0$
\begin{eqnarray*}
I_1&=&\reip\{\Pi_{1}((r,\infty)\times\{ e_{1}\} )\neq 0, \Pi_{2}((tr,\infty)\times\{e_{2}\})=0\}\\
&\leq&\reip\{\xi_{\alpha,\sigma}\in (r,\infty)\times B\}\\
&\leq&\reip\{\pi_{1}((r,\infty)\times\{e_{1}\})\neq 0\}\\
&=&I_2,
\end{eqnarray*}
où \[I_1=(1-\exp(-pr^{-\alpha}))\exp(-q(tr)^{-\alpha}) \;\; \mbox{et} \;\; I_2=1-\exp(-pr^{-\alpha}).\]
Les deux termes $I_1$ et $I_2$ sont équivalents à $pr^{-\alpha}$ lorsque $r\rightarrow\infty$, ainsi on a 
\[\reip\{\xi_{\alpha,\sigma}\in B\times (r,\infty)\}\sim pr^{-\alpha} \;\;\;  \mbox{lorsque} \;\;\;  r\rightarrow\infty.\]  
De même manière, pour $B=[(s,1), (t,1)], 0\leq s<b<t\leq 1$ on a $\reip\{\xi_{\alpha,\sigma}\in (r,\infty)\times B\}\sim qr^{-\alpha} \  \mbox{si} \  r\rightarrow\infty$. 
La loi de $\xi_{\alpha,\sigma}$ vérifie la condition de variation régulière (\ref{regulier2}).


\vspace{0.5cm}
\subsection{Lois $\alpha$-stables}\label{subsect-parastable}
La difficulté dans l'étude de lois $\alpha$-stables est que, sauf pour les cas particuliers ($\alpha=0.5,1 \, \mbox{et} \, 2$), il n'y a pas de forme explicite pour la densité. Les vecteurs aléatoires $\alpha$-stables sont déterminés par leur fonction caractéristique.

\begin{definition}\label{deffc}
La fonction caractéristique d'un vecteur aléatoire $\alpha$-stable ($0<\alpha<2$) dans $\rd$ s'exprime par l'expression suivante:
\begin{equation*}
\phi_{\alpha,\sigma}(t)=\exp\left(-\frac{1}{C_{\alpha}}\ints\psi_\alpha(\langle t,s \rangle)\sigma(ds)+i\langle \delta, t\rangle \right),
\end{equation*}
où $\sd=\{x\in\rd \, | \, \|x\|=1\}$, $\sigma$ est une mesure finie sur la sphère unité $\sd$, $\delta$ est un vecteur dans $\rd$, $\langle t,s \rangle$ représente le produit scalaire,
\begin{equation}\label{ca}
C_{\alpha}=\left\{\begin{array}{cr}
\frac{1-\alpha}{\Gamma(2-\alpha)\cos(\pi\alpha/2)}, &\alpha\neq 1,\\
\frac{2}{\pi}, & \alpha=1,
\end{array}\right.
\end{equation}
et
\[\psi_\alpha(x)=\left\{\begin{array}{lr}
 |x|^\alpha\left(1-i\,\mbox{sign}(x)\tan\displaystyle\frac{\pi\alpha}{2}\right),&\alpha\neq 1,\\
|x|\left(1+i \displaystyle\frac{\pi}{2}\,\mbox{sign}(x)\ln|x|\right),& \alpha = 1.
\end{array}\right.\]
\end{definition}
La loi stable dans $\rd$ est spécifiée par un nombre $\alpha$ entre $0$ et $2$ dit {\it indice de stabilité}, une mesure finie $\sigma$ sur $\sd$ dite {\it mesure spectrale} et un vecteur $\delta$ dans $\rd$. Il existe des méthodes numériques pour obtenir la densité et la fonction de répartition approximatives d'une loi stable \cite{Nolan97}. 

Dans le cas unidimensionnel la sphère unité ne contient que deux points, i.e. $S^0=\{-1,1\}$. La mesure spectrale $\sigma$ se réduit à deux valeurs $\sigma(-1)$ et $\sigma(1)$. La fonction caractéristique peut être écrite comme suit
\[\phi_{\alpha,\sigma}(t)=\left\{\begin{array}{lr}
 \exp\left(-\frac{\gamma}{C_{\alpha}} |t|^\alpha\left(1-i\beta\,\mbox{sign}(t)\tan\displaystyle\frac{\pi\alpha}{2}\right)+i\delta t\right),&\alpha\neq 1,\\
\exp\left(-\frac{\gamma}{C_{\alpha}} |t|\left(1+i\beta \displaystyle\frac{\pi}{2}\,\mbox{sign}(t)\ln|t|\right)+i\delta t\right),& \alpha = 1,
\end{array}\right.\]
où $\gamma=\sigma(1)+\sigma(-1)$ et $\beta=(\sigma(1)-\sigma(-1))/\gamma$. La loi stable unidimensionnelle est déterminée par quatre paramètres: $\alpha\in(0,2]$, $\beta\in[-1,1]$, $\gamma>0$ et $\delta\in \Bbb{R}$. Le paramètre $\beta$ contrôle l'asymétrie. La loi est dite {\em biaisée totalement vers la droite} si $\beta=1$ et {\em biaisée totalement vers la gauche} si $\beta=-1$. Si $\beta=0$ la densité est symétrique par rapport à $\delta$. Les paramètres $\gamma$ et $\delta$ sont les paramètres d'échelle et de position. Nous utilisons la notation $\mathcal{S}_{d}(\alpha,\sigma,\delta)$ pour la loi $\alpha$-stable $d$-dimensionnelle avec $d\geq 2$, la notation $\mathcal{S}_{1}(\alpha, \beta,\gamma,\delta)$ pour $d=1$.

La paramétrisation qu'on rencontre souvent dans les applications est celle de Samorodnitsky et Taqqu \cite{Samorodnitsky94}. La différence entre notre paramétrisation et celle de \cite{Samorodnitsky94} concerne la masse totale de la mesure spectrale $\sigma$. Notons $\sigma_{s}$ le paramètre d'échelle si $d=1$ et la mesure spectrale si $d\geq 2$ dans \cite{Samorodnitsky94}. Dans le cas unidimensionnel
\[\gamma=C_{\alpha}\sigma_{s}^{\alpha},\]
et dans le cas multidimensionnel
\[\sigma(\cdot)=C_{\alpha}\sigma_{s}(\cdot),\]
où $C_{\alpha}$ est définie par (\ref{ca}). L'avantage de notre paramétrisation est que la mesure spectrale coïncide avec celle qui apparaît dans la propriété de régularité (\ref{regulier2}). Soit $X$ la variable aléatoire de loi $\mathcal{S}_{1}(\alpha, \beta,\gamma,\delta)$, alors on a
\begin{equation}\label{powerlaw1d}
\left\{\begin{array}{lll}\lim\limits_{x\rightarrow\infty}x^\alpha\reip\{X>x\}&=\sigma(1),\\
\lim\limits_{x\rightarrow\infty}x^\alpha\reip\{X<-x\}&=\sigma(-1),\end{array}\right.
\end{equation}
où
\[\sigma(1)=\frac{(1+\beta)}{2}\gamma, \;\; \sigma(-1)=\frac{(1-\beta)}{2}\gamma.\] 
La démonstration de (\ref{powerlaw1d})  peut être trouvée dans \cite{Samorodnitsky94}, page 16 et 197.

\begin{description}
\item[Propriété 1.] Soit $X$ le vecteur aléatoire de loi $\mathcal{S}_{d}(\alpha,\sigma,\delta)$. Le vecteur aléatoire $X$ est \sas \, dans $\rd$ avec $0<\alpha\leq 2$ si et seulement si
\begin{description}
\item[(a)] $\alpha\neq 1$, $\delta=0$,
\item[(b)] $\alpha=1$, $\int_{\sd}s^{(k)}\sigma(ds)=0$, pour $k=1,2,\ldots,d$.
\end{description}

\item[Propriété 2.]Le vecteur aléatoire $\alpha$-stable est {\em symétrique} si et seulement si $\delta=0$ et $\sigma$ est une mesure symétrique sur $\sd$, i.e. 
\begin{equation}\label{symsig}
\sigma(A)=\sigma(-A)
\end{equation}
pour tout ensemble borélien $A$ dans $\sd$. 
\end{description}
\vspace{0.5cm}

\subsubsection*{Données $\alpha$-stables simulées}
Il existe des méthodes connues pour simuler les vecteurs aléatoires stables. Chambers et al. \cite{Chambers76} ont proposé une méthode pour simuler les variables aléatoires stables unidimensionnelles arbitraires. Plus récemment une technique précise a été construite dans \cite{Weron96} en utilisant une transformation non-linéaire d'une paire des variables aléatoires indépendantes de loi uniforme et de loi exponentielle. Une méthode pour simuler les vecteurs aléatoires stables basée sur la discrétisation de la mesure spectrale a été présentée par Modarres et Nolan \cite{Modarres94}. Le résultat suivant présenté par R. LePage et al. montre qu'on peut simuler des lois stables en utilisant la représentation de série.

\begin{theoreme}\label{lepagerd} (\cite{Lepage81}, Th. 3) 
La série définie par (\ref{lepageintro}) avec $\alpha<1$ ou $\sigma$ symétrique  a la loi \sas.
\end{theoreme}

Soit $\xi_{\alpha,\sigma}$ un vecteur aléatoire stable d'indice $\alpha$ et de mesure spectrale $\sigma$ vérifiant la condition du théorème précédent. On simule $\xi_{\alpha,\sigma}$ en utilisant la somme partielle de la série de LePage (\ref{lepageintro}),
\begin{equation}\label{somparstab}
\hat\xi_{k}=\sum_{i=1}^k c\Gamma_i^{-1/\alpha}\epsilon_i.
\end{equation}

\begin{figure}[!htbp]
\hspace{-0.5cm}\begin{tabular}{c}
\includegraphics[width=80mm, height=61mm]{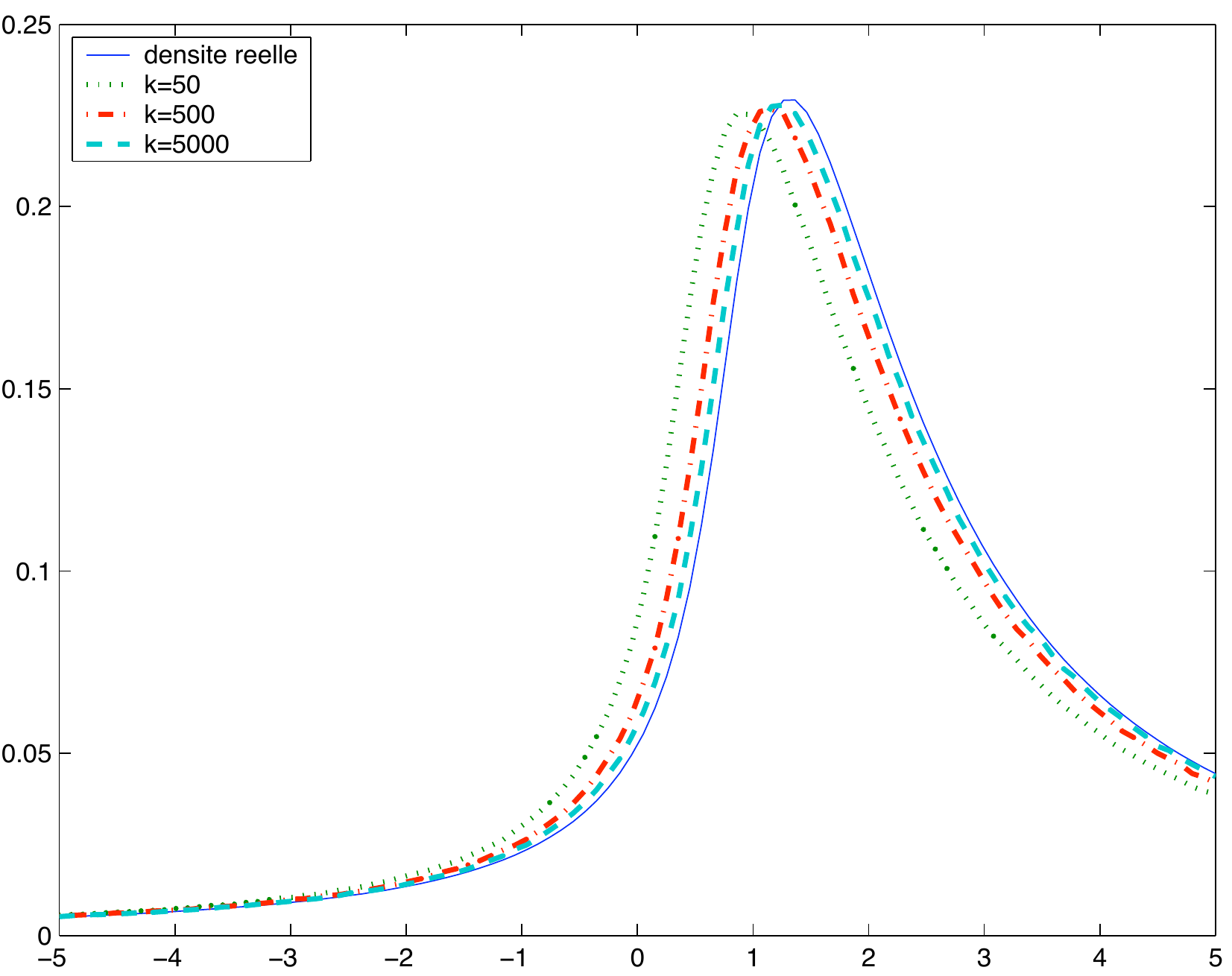} \\
\includegraphics[width=80mm, height=60mm]{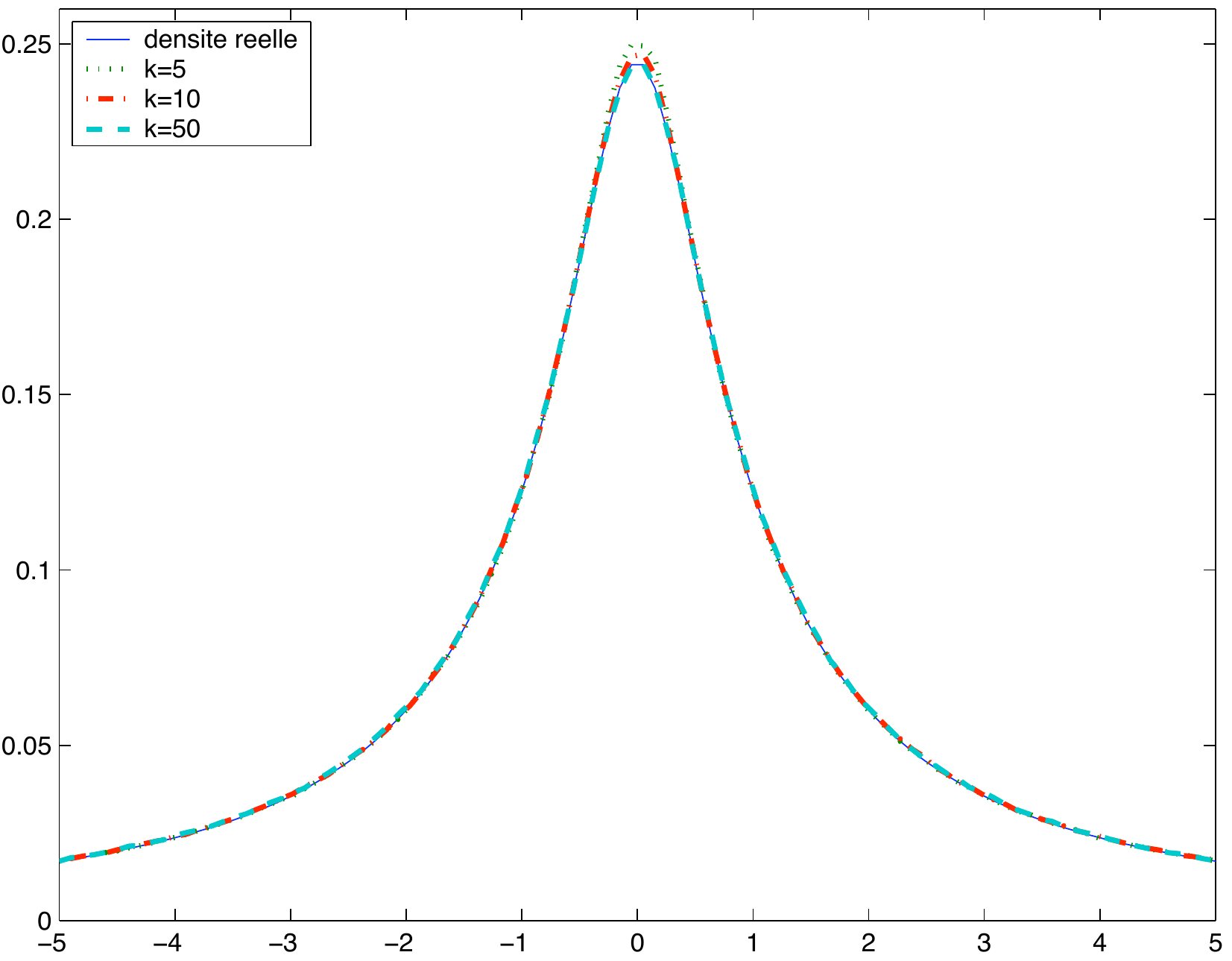} 
\end{tabular}
\caption[Convergence de la série de LePage.]{Densités des lois
  $\mathcal{S}_{1}(0.75,0.5, 1,0)$ (haut) et
  $\mathcal{S}_{1}(0.75,0, 1,0)$ (bas) simulées par la somme
  partielle de série de LePage avec différentes valeurs de $k$. Le nombre de données simulées est $10^7$.}
\label{densimustab1d}
\end{figure}

Dans un premier temps deux lois unidimensionnelles $\mathcal{S}_{1}(0.75,0.5,1,0)$ (non symétrique) et $\mathcal{S}_{1}(0.75,0, 1,0)$ (symétrique) sont considérées. On simule les variables aléatoires en utilisant la somme (\ref{somparstab}) avec le nombre de termes $k$ différent. Les valeurs $50$, $500$ et $5000$ sont utilisées pour simuler la loi stable non symétrique, tandis que la loi stable symétrique est simulée avec $k= 5, 10$ et $50$. On observe que l'approche de la densité empirique à la densité réelle est plus rapide pour la loi stable symétrique, voir la figure \ref{densimustab1d}. Cela signifie que la somme partielle (\ref{somparstab}) donne une bonne approximation pour les lois $\alpha$-stables symétriques. Dans le cas non-symétrique, la convergence est lente pour une valeur de $k$ modeste. 

\begin{figure}[!htbp]
\begin{center}
\begin{tabular}{cc}
a.\includegraphics[width=50mm, height=50mm]{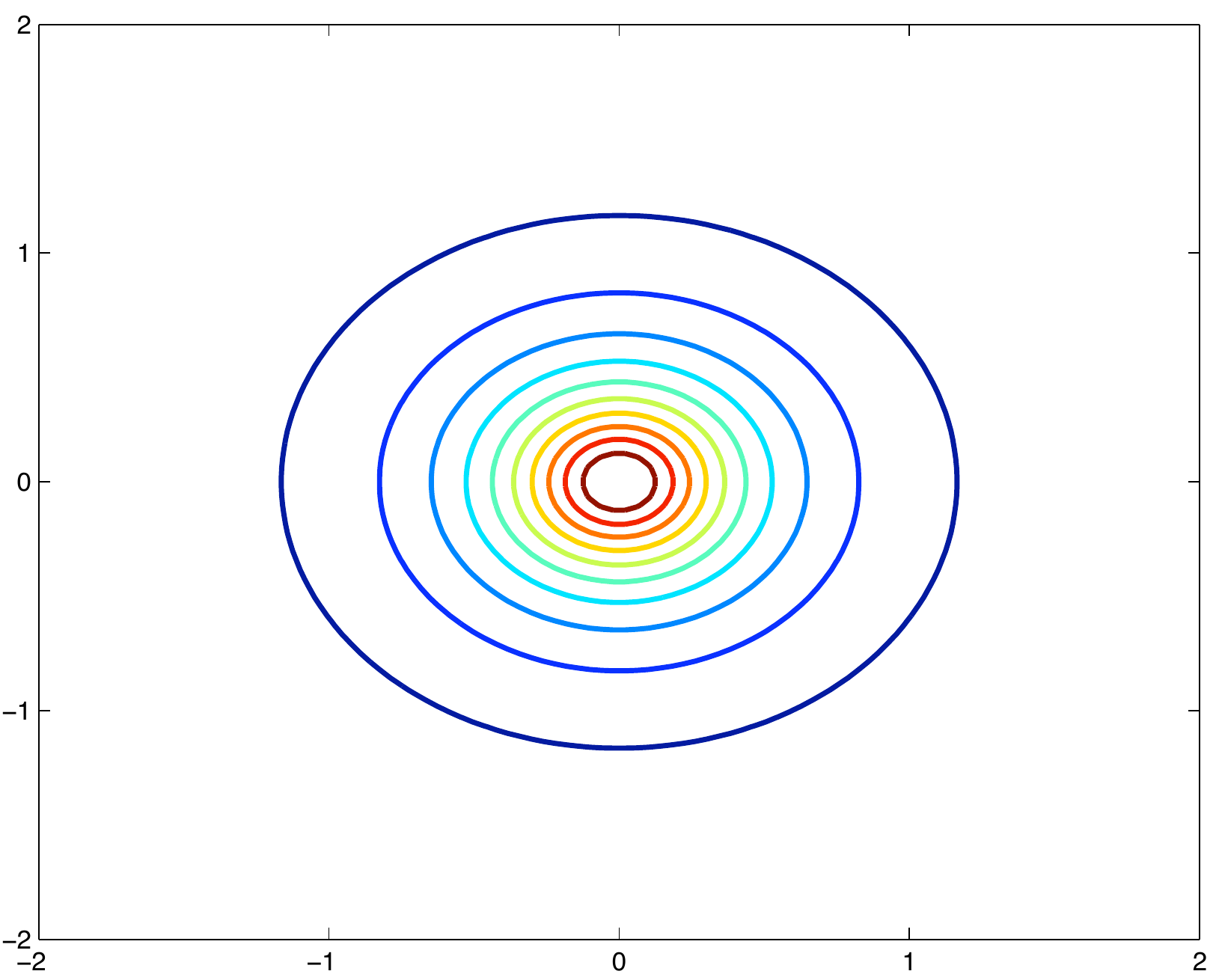} & 
 \includegraphics[width=50mm, height=50mm]{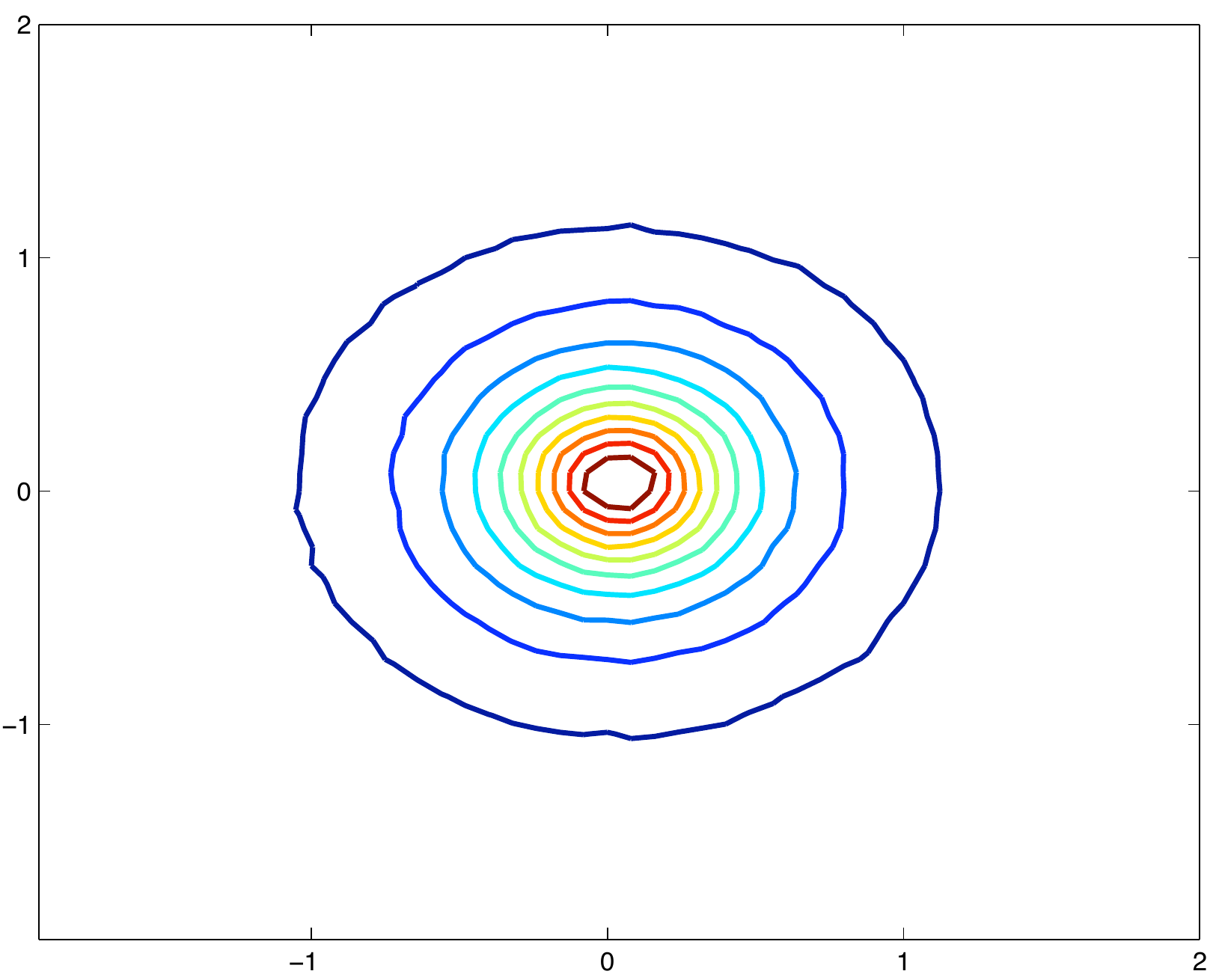} \\
b. \includegraphics[width=50mm, height=50mm]{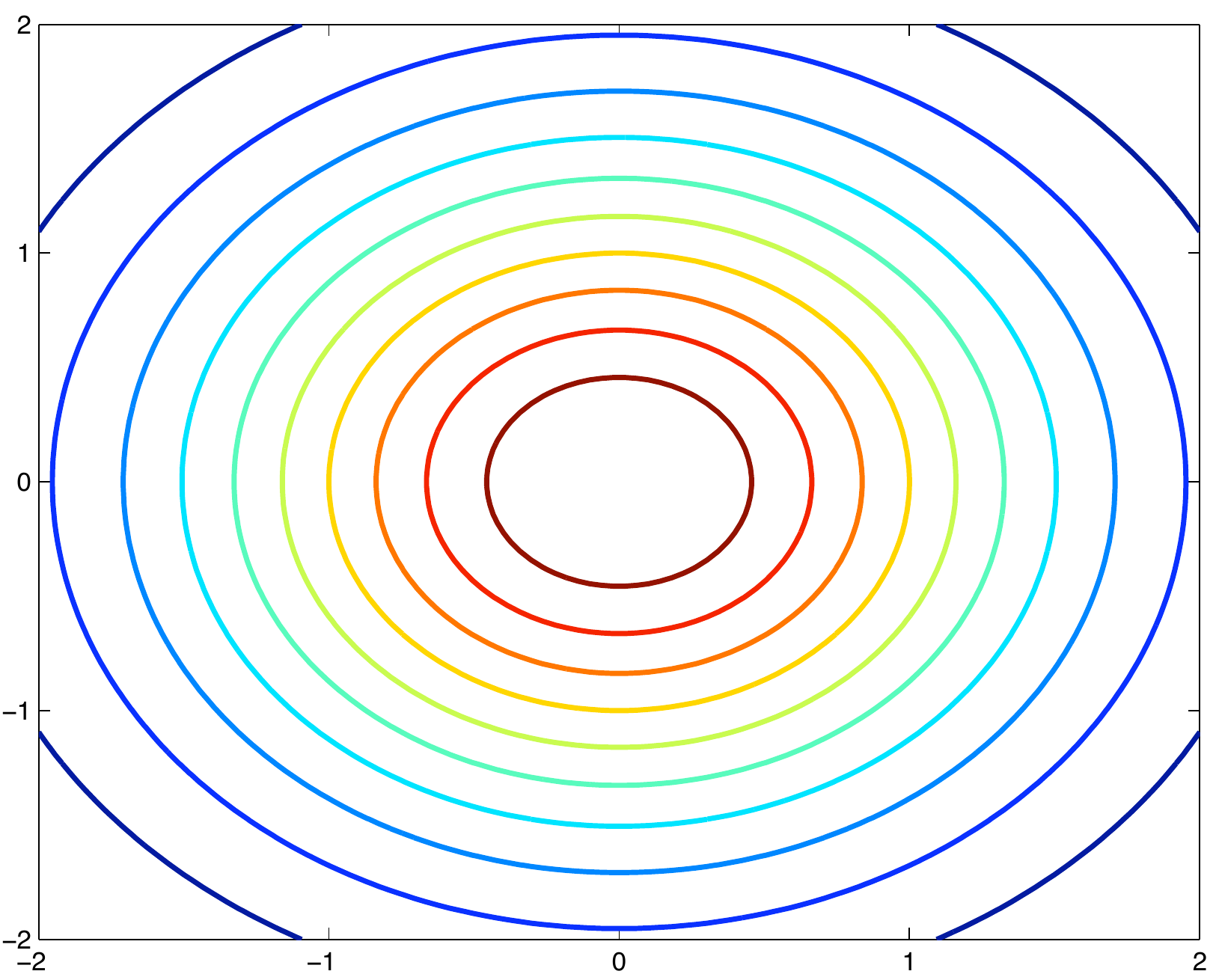} & 
 \includegraphics[width=50mm, height=50mm]{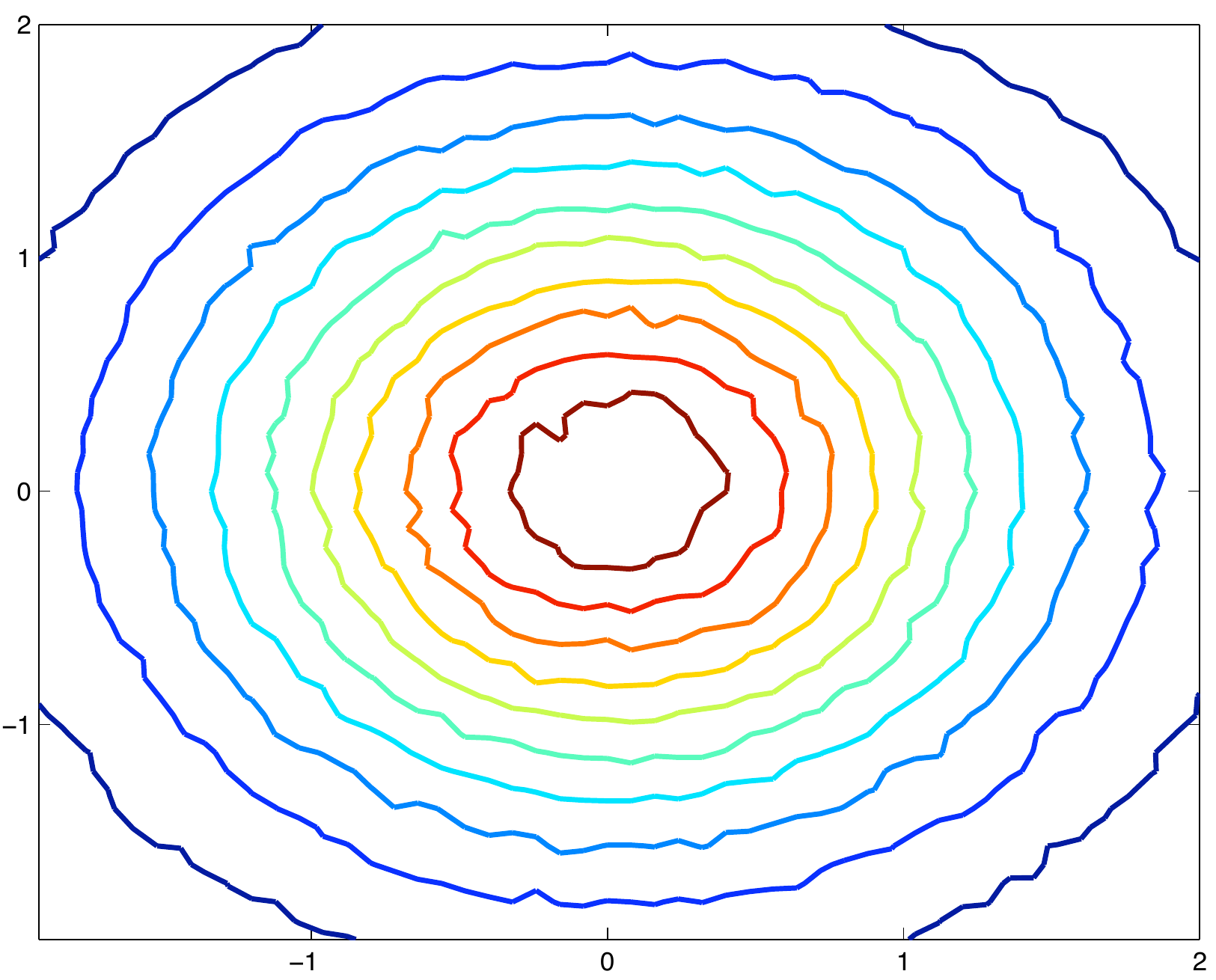}
\end{tabular}
\end{center}
\caption[Simulations des lois $\alpha$-stables isotropiques bivariées.]{Courbes de niveau des densités réelles (gauche) et empiriques (droite) de lois $\alpha$-stables bivariées de mesure spectrale uniforme, $\alpha=0.75$ (a.), $1.5$ (b.), $k=10$. Le nombre de données simulées est $10^7$.}
\label{densimustab2diso}
\end{figure}

Ensuite les vecteurs aléatoires $2$-dimensionnels de lois $0.75$ et $1.5$-stables avec la mesure spectrale uniforme sont simulés. Le nombre de termes de la somme partielle est $10$. On calcule la densité réelle en utilisant la fonction Matlab ``mvstablepdf'' dans le package appelé ``STABLE'' obtenu sur demande au site \url{http://www.RobustAnalysis.com/}. Les courbes de niveau des densités réelles et empiriques sont présentées dans la figure \ref{densimustab2diso}. 

\begin{figure}[!htbp]
\begin{center}
\begin{tabular}{cc}
a. \includegraphics[width=50mm, height=50mm]{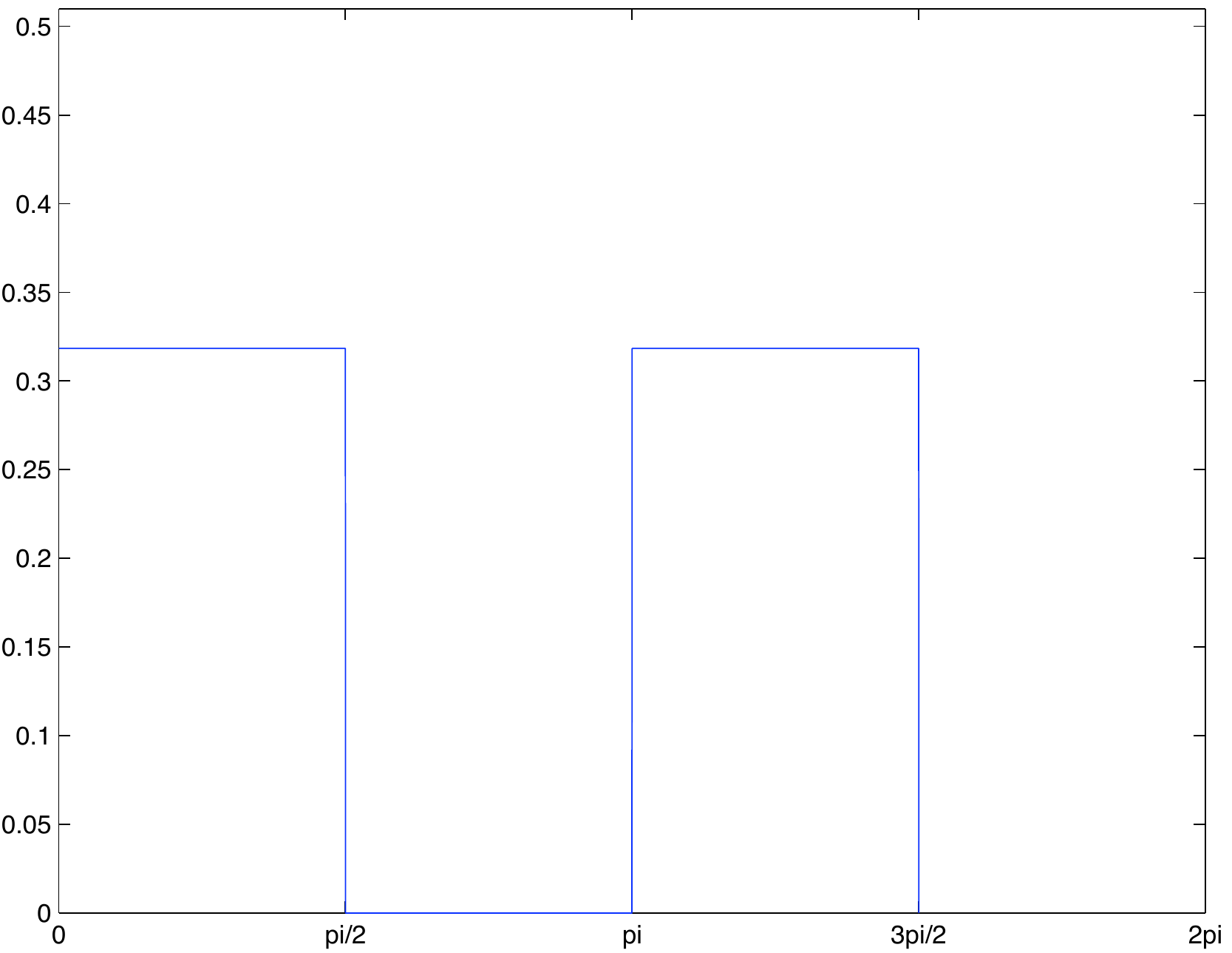} & 
\includegraphics[width=50mm, height=50mm]{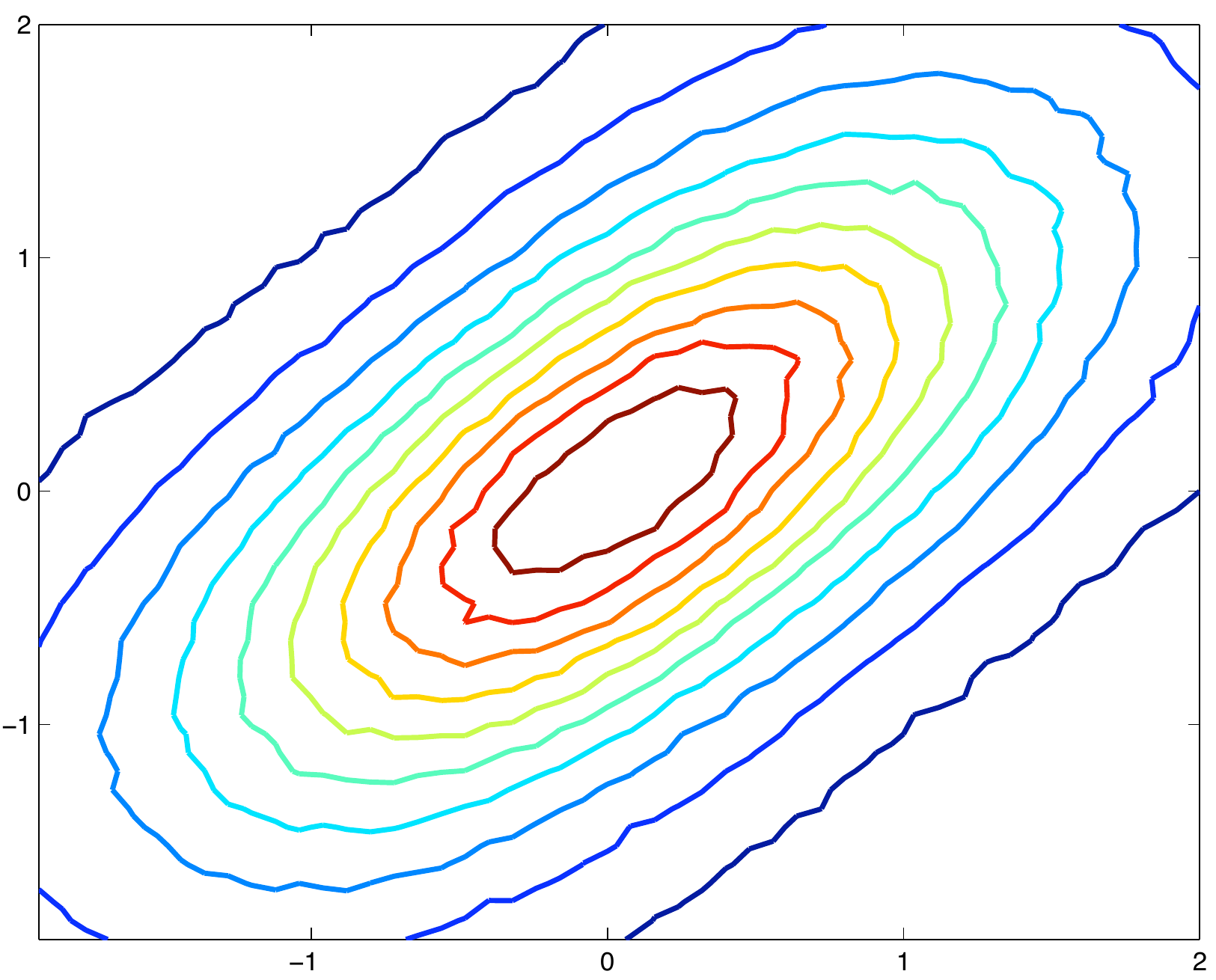} \\
b. \includegraphics[width=50mm, height=50mm]{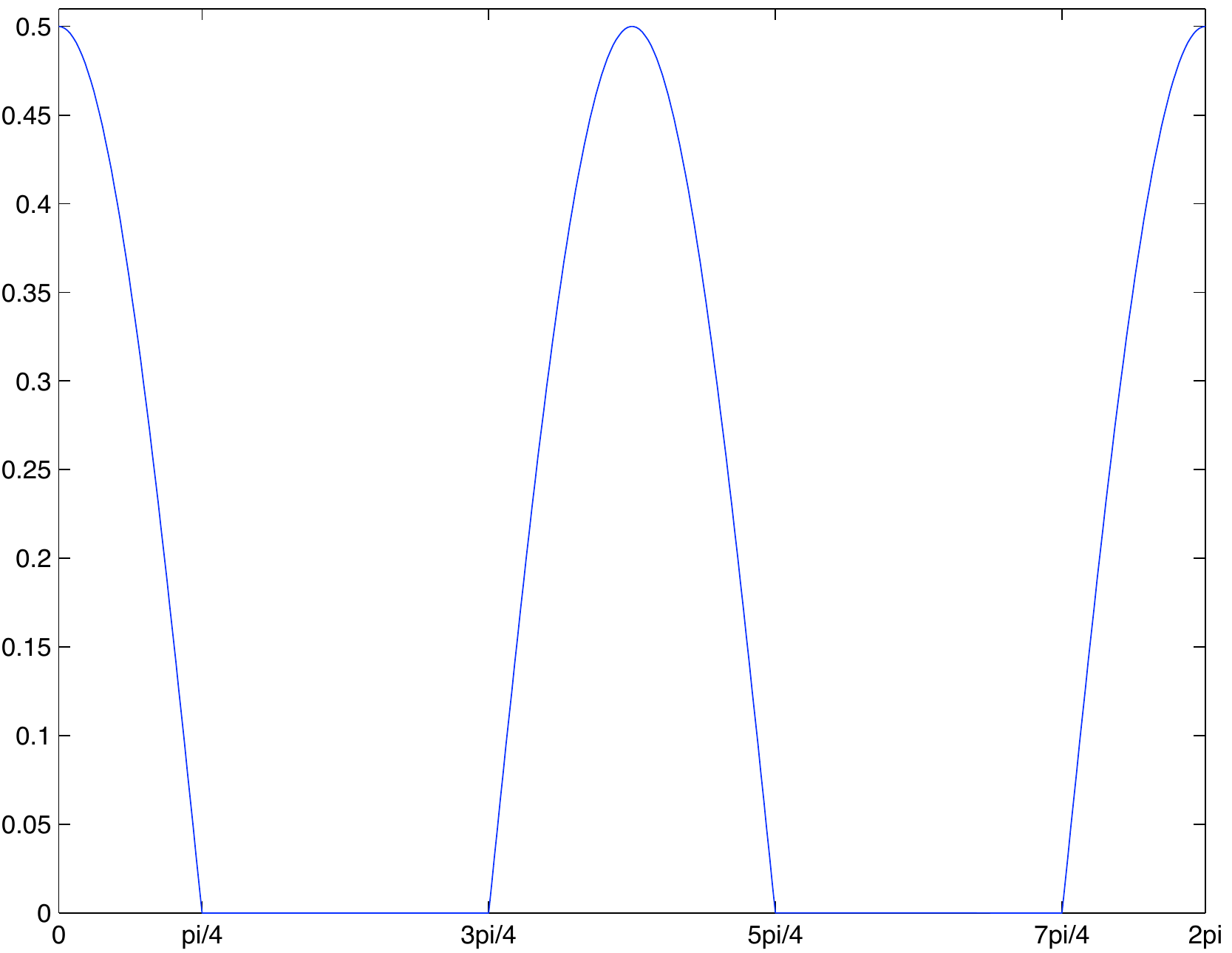} &
\includegraphics[width=50mm, height=50mm]{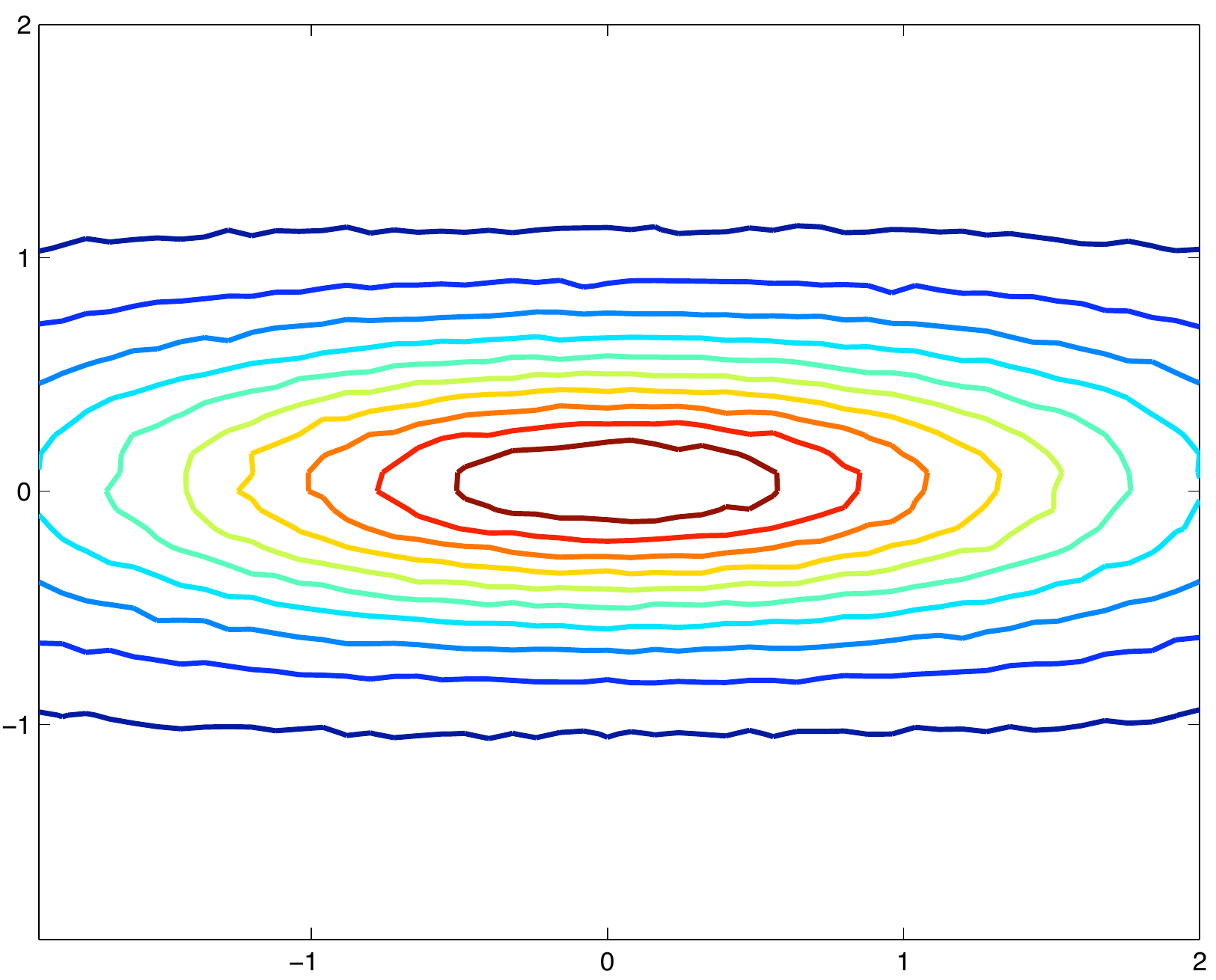} \\
c. \includegraphics[width=50mm, height=50mm]{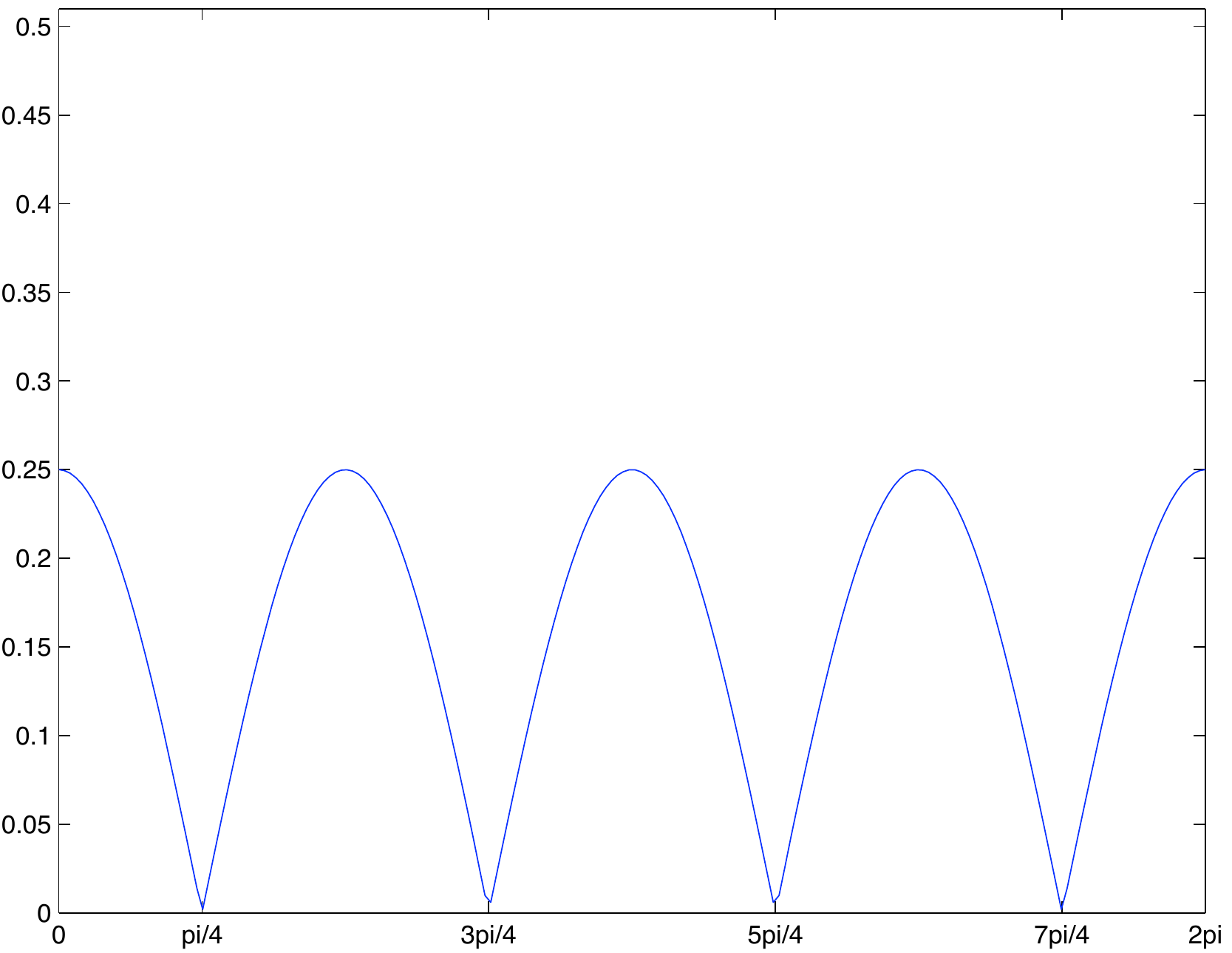} &
\includegraphics[width=50mm, height=50mm]{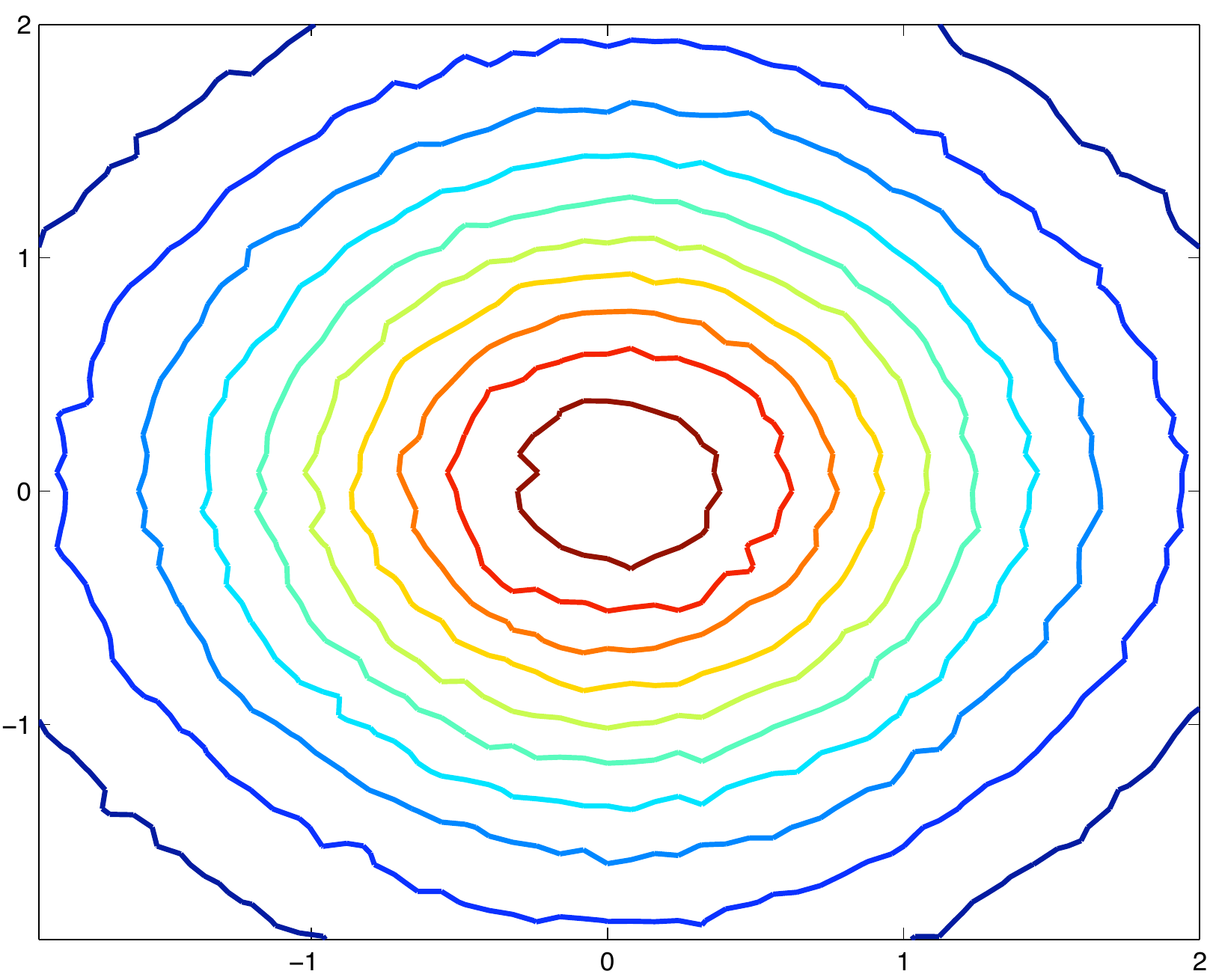}
\end{tabular}
\end{center}
\caption[Simulations des lois $\mathcal{S}_{2}(1.5,\sigma_{i},0),  i=1,2,3$.]{Densités de la mesure spectrale $\sigma_{i},\,  i=1,2,3$ définies par (\ref{defsigma1}), (\ref{defsigma2}) et (\ref{defsigma3}) (gauche) et courbes de niveau des densités empiriques de lois $\mathcal{S}_{2}(1.5,\sigma_{i},0),  i=1,2,3$ (droite), $k=10$. Le nombre de données simulées est $10^7$. }
\label{densimustab2d}
\end{figure}

Enfin, on simule les vecteurs aléatoires $2$-dimensionnels de lois $1.5$-stables avec les mesures spectrales $\sigma_{1}$, $\sigma_{2}$ et $\sigma_{3}$ dont les densités sont définies par 
\begin{eqnarray}
f_{1}(\theta)&=&\left\{\begin{array}{ll}
\frac{1}{\pi}, & \theta\in[0, \frac{\pi}{2})\cup[\pi, \frac{3\pi}{2}),\\
0,& \mbox{sinon},
\end{array}\right.\label{defsigma1}\\
f_{2}(\theta)&=&\left\{\begin{array}{ll}
\frac{1}{2}\cos(2\theta), & \theta\in[-\frac{\pi}{4}, \frac{\pi}{4})\cup[\frac{3\pi}{4}, \frac{5\pi}{4}),\\
0, & \mbox{sinon},
\end{array}\right.\label{defsigma2}\\
f_{3}(\theta)&=&\frac{1}{4}|\cos (2\theta)|,\;\; \theta\in [0,2\pi).\label{defsigma3}
\end{eqnarray}
Ces trois mesures vérifient la condition (\ref{symsig}), les lois strictement stables correspondantes sont donc symétriques. Le package de programmes qu'on utilise ne permet pas de calculer la densité de la loi stable avec la mesure spectrale absolument continue. Nous présentons dans la figure \ref{densimustab2d}  les densités de mesure spectrale $f_i$ et les courbes de niveau des densités empiriques de lois $\mathcal{S}_{2}(1.5,\sigma_{i},0), i=1,2,3$.

Un des usages de ces données est de tester la robustesse des procédures statistiques multivariées: engendrer aléatoirement l'ensemble de données et évaluer les statistiques par la procédure considérée. Un autre usage possible est de calculer la probabilité $\reip \{X\in A\}$ pour l'ensemble $A\subset \Bbb{R}^d$ où $X$ est le vecteur aléatoire qu'on simule. En général le calcul numérique de cette probabilité est difficile, voir par exemple \cite{Nolan95} et \cite{Nolan97}. On peut les estimer par la méthode standard de Monte-Carlo une fois que nous savons générer les vecteurs aléatoires de loi prescrite. Le calcul devient coûteux si on prend le nombre de termes $k$ grand. Donc cette méthode n'est pas avantageuse pour la simulation des lois stables non symétriques. Le troisième usage de ces vecteurs est la simulation dans les domaines d'application, par exemple l'analyse du portefeuille stable \cite{Press72}.

\section{Conclusion}\label{conclusion}

Sous l'hypothèse de variation régulière tous les éléments \sas s avec $\alpha\in (0,1)$ peuvent être représentés par la somme des points d'un processus ponctuel poissonnien, dite représentation de LePage. La représentation de LePage de la loi strictement max-stable implique que la loi ne dépend que des points les plus grands du processus ponctuel poissonnien correspondant. En utilisant cette propriété, on peut déduire les lois jointes et marginales des vecteurs aléatoires max-stables avec les mesures spectrales données. La propriété de variation régulière de la queue des lois est vérifiée. Les éléments aléatoires \sas s vérifiant la condition de variation régulière peuvent être simulés en utilisant la somme partielle de la série de LePage. Le temps de calcul de simulation est faible pour les lois strictement max-stables et \sas s symétriques. Les résultats de simulation montrent visuellement que l'indice $\alpha$ et la mesure spectrale $\sigma$ mesurent respectivement la disposition des masses et la disposition directionnelle de probabilité autour du centre.    

\bibliographystyle{plain}
\bibliography{reference.bib}

\end{document}